\documentclass[aop,MSNbibl,citesort,dvips]{arximspdf}
\usepackage{graphicx}
\doi{10.1214/10-AOP639}
\volume{40}
\issue{3}
\pubyear{2012}
\firstpage{1069}
\lastpage{1104}

\makeatletter

\newcommand{\mes}{\operatorname{mes}}

\newcommand{\cA}{\mathcal A}
\newcommand{\cB}{\mathcal B}
\newcommand{\cC}{\mathcal C}
\newcommand{\cD}{\mathcal D}
\newcommand{\kk}{\mathbf k}
\newcommand{\zz}{\mathbf z}

\newtheorem{theorem}{Theorem}
\newtheorem{corollary}{Corollary}
\newtheorem{lemma}{Lemma}
\newtheorem{proposition}{Proposition}
\newproclaim{remark}{Remark}

\makeatother

\begin{document}
\begin{frontmatter}

\title{Tail approximations of integrals of Gaussian random fields}
\runtitle{Tail approximations of integrals of Gaussian random fields}

\begin{aug}
\author[A]{\fnms{Jingchen} \snm{Liu}\corref{}\thanksref{t1}\ead[label=e1]{jcliu@stat.columbia.edu}}
\runauthor{J. Liu}
\affiliation{Columbia University}
\address[A]{Columbia University\\
1255 Amsterdam Ave.\\
New York, New York 10027\\
USA\\
\printead{e1}} 
\end{aug}

\thankstext{t1}{Supported in
part by the Institute of Education Sciences, U.S. Department of
Education, through Grant R305D100017 and by NSF CMMI-1069064.}

\received{\smonth{3} \syear{2010}}
\revised{\smonth{12} \syear{2010}}

%
\begin{abstract}
This paper develops asymptotic approximations of $P(\int_T e^{f(t)}\,dt >
b)$ as $b\rightarrow\infty$ for a homogeneous smooth Gaussian random
field, $f$, living on a compact $d$-dimensional Jordan measurable set
$T$. The integral of an exponent of a Gaussian random field is an
important random variable for many generic models in spatial point
processes, portfolio risk analysis, asset pricing and so forth.

The analysis technique consists of two steps: 1. evaluate the tail
probability $P(\int_\Xi e^{f(t)}\,dt > b)$ over a small domain $\Xi$
depending on $b$, where $\mes(\Xi) \rightarrow0$ as $b\rightarrow
\infty$ and $\mes(\cdot)$ is the Lebesgue measure; 2. with $\Xi$
appropriately chosen, we show that $P(\int_T e^{f(t)}\,dt > b) =(1+o(1))
\mes(T)\times \mes^{-1}(\Xi) P(\int_\Xi e^{f(t)}\,dt > b)$.
\end{abstract}

%

\end{frontmatter}

\section{Introduction}\label{SecIntro}
We consider a Gaussian random field living on a $d$-di-\break mensional
domain $T\subset R^d$, $\{f(t)\dvtx t\in T\}$. For every
finite subset $\{t_1,\ldots,\allowbreak t_n\}\subset T$, $(f(t_1),\ldots,f(t_n))$
is a multivariate Gaussian random vector. The quantity of
interest is
\[
P\biggl(\int_T e^{f(t)}\,dt
> b \biggr)
\]
as $b\rightarrow\infty$.

The motivations of the study of $\int_T e^{f(t)}\,dt$ are from
multiple sources. We will
present a few of them. Consider a point process on $R^d$
associated with a~Poisson random measure $\{N_A\}_{A\in
\mathcal B}$ with intensity $\lambda(t)$, where $\mathcal B$
represents the Borel sets of $R^d$. One important task in spatial
modeling is to build in dependence structures. A~popular
strategy is to let $f(t)= \log\lambda(t)$, which can take all
values in $R$, and model $f(t)$ as a Gaussian random field.
Then,\vadjust{\goodbreak} $\int_A e^{f(t)}\,dt = E(N(A)|\lambda(\cdot))$ for all $A\in
\mathcal B$. With the multivariate Gaussian structure, it is easy to
include linear predictors in the intensity process. For instance,
\cite{ChLe95} models $f(t) = U(t) + W(t)$, where $U(t)$ is the observed
(deterministic) covariate process and $W(t)$ is a stationary AR(1)
process. Similar models can be found in \cite{DDW00,Camp94,Zeger88}
which are special cases of the Cox process \cite{COX55,COIS80}. Such
a modeling approach has been applied to many disciplines, a~short list
of which is as follows: astronomy, epidemiology, geography, ecology,
material science and so forth.

In portfolio risk analysis, consider a portfolio consisting of
equally weighted assets $(S_1,\ldots,S_n)$. One stylized model is that
$(\log S_1,\ldots,\log S_n)$ is a multivariate normal random vector
(cf. \cite{DufPan97,Ahs78,BasSha01,GHS00,Due04}). The value of the
portfolio $S= \sum_{i=1}^n S_i$ is then the sum of correlated
log-normal random variables. If
one can represent each asset price by the value of a Gaussian
random field at one location in $T$, that is, $\log S_i
=f(t_i)$. As the portfolio size tends to infinity and the asset
prices become more correlated, the limit of the unit share
price of the portfolio is $\lim_{n\rightarrow\infty} S/n =
\int_T e^{f(t)}\,dt$. For more general cases, such as unequally
weighted portfolios, the integral is possibly with respect to
some other measures instead of the Lebesgue measure.

In option pricing, if we let $S(t)$ be a geometric Brownian motion (cf.
\cite{BlaSch73}, Chapter 5 of \cite{Duf01}, Chapter 3.2 of \cite
{fGLA04a}), the payoff function of an Asian option (with expiration
time $T$) is a function of $\int_0^T S(t) \,dt$. For instance, the
payoff of an Asian
call option with strike price $K$ is $\max(\int_0^T S(t)
\,dt-K,0)$; the payoff of a digital Asian call option is
$I(\int_0^T S(t) \,dt>K)$.

We want to emphasize that the extreme behavior of $\int_T
e^{f(t)}\,dt$ connects closely to that of $\sup_T f(t)$. As we
will show in Theorem \ref{ThmM}, with the threshold~$u$
appropriately chosen according to $b$, the probabilities of
events\break $\{\int_T e^{f(t)}\,dt>b \}$ and $\{\sup_T f(t)\,dt>u \}$
have asymptotically the same decaying rate. It suggests that
these two events have substantial overlap with each other.
Therefore, we will borrow the intuitions and existing
results on the high excursion of the supremum of random
fields for the analysis of $\int_T e^{f(t)}\,dt$.

There is a vast literature on the extremes of Gaussian random
fields mostly focusing on the tail probabilities of $\sup_T f(t)$ and
its associated geometry. The results contain general bounds on $P(\sup
_T f(t) >b)$ as well as sharp asymptotic approximations as
$b\rightarrow\infty$. A partial literature contains
\cite{LS70,MS70,ST74,Bor75,Bor03,LT91,TA96,Berman85}. Several
methods have been introduced to obtain asymptotic
approximations, each of which imposes different regularity
conditions on the random fields. A~few examples are given as
follows. The double sum method \cite{Pit95} requires
expansions of the covariance function and locally stationary structure. The
Euler--Poincar\'{e} Characteristic of the excursion set
[$\chi(A_b)$] approximation uses the fact that $P(M>b)\approx
E(\chi(A_b))$, which requires the random field to be at least
twice differentiable \cite{Adl81,TTA05,AdlTay07}. The tube
method \cite{Sun93} uses the Karhunen--Lo\`{e}ve expansion and
imposes differentiability\vadjust{\goodbreak} assumptions on the covariance
function (fast decaying eigenvalues). The Rice method \cite{AW05,AW08,AW09}
represents the distribution of $M$ (density function) in an
implicit form. Recently, the efficient simulation algorithms
are explored by \cite{ABL08,ABL09}. These two papers provided
computation schemes that run in polynomial time to compute the
tail probabilities for all H\"{o}lder continuous Gaussian random
fields and in constant time for twice differentiable and
homogeneous fields. In addition, \cite{AST09} studied the geometric
properties of a high level
excursion set for infinitely divisible non-Gaussian fields as well as the
conditional distributions of such properties given the high excursion.

The distribution of $\int e^{f(t)}\,dt$ for the special case that $f(t)$
is a Wiener process has been studied by \cite{Yor92,Duf01}. For other
general functionals of Gaussian processes and multivariate Gaussian
random vectors, the tail approximation of the finite sum of correlated
log-normal
random variables has been studied by \cite{AR08}. The
corresponding simulation is studied in \cite{BJR10}. The gap
between the finite sums of log-normal r.v.'s and the integral
of continuous fields is substantial in the aspects of both
generality and techniques.

The basic strategy of the analysis consists of two steps. The
first step is to partition the domain $T$ into $n$ small
squares of equal size denoted by $A_i$, $i=1,\ldots,n$, and
develop asymptotic approximations for each $p_i =P(\int_{A_i}
e^{f(t)}\,dt > b)$. The size of $A_i$ will be chosen carefully
such that it is valid to use Taylor's expansion on $f(t)$ to
develop the asymptotic approximations of $p_i$. The second
step is to show that $P(\int_T e^{f(t)}\,dt>b)= (1+o(1))
\sum_{i=1}^n p_i$. This implies that when computing $P(\int_T
e^{f(t)}\,dt>b)$, we can pretend that all the $\int_{A_i}
e^{f(t)}\,dt$'s are independent, though they are truly highly
dependent. The sizes of the $A_i$'s need to be chosen
carefully. If $A_i$ is too large, Taylor's expansion may not be
accurate; if $A_i$ is too small, the dependence of the fields
in different $A_i$'s will be high and the second step
approximation may not be true. Since the first step of the
analysis requires Taylor's expansion of the field, we will need
to impose certain conditions on the field, which will be given
in Section \ref{SecPre}.

This paper is organized as follows. In Section \ref{SecPre} we
provide necessary background and the technical conditions on the
Gaussian random field in context. The main theorem and its
connection to asymptotic approximation of $P(\sup_T f(t)
>b)$ are presented in Section \ref{SecMain}. In addition, two
important steps of the proof are given in the same section, which lay
out the proof strategy. Sections \ref{sec_proof1} and
\ref{sec_proof2} give the proofs of the two steps
presented in Section \ref{SecMain}. Detailed lemmas and their proofs
are given in the \hyperref[app]{Appendix}.

\section{Some useful existing results}\label{SecPre}

\subsection{Preliminaries and technical conditions for Gaussian random field}

Consider a homogeneous Gaussian random field, $f(t)$, living on
a domain $T$. Denote the covariance function by
\[
C(t-s) =
\operatorname{Cov}(f(s), f(t)).
\]
Throughout this paper, we assume that the
random field satisfies the following conditions:
\begin{longlist}[(C2)]
\item[(C1)] $f$ is homogeneous with $Ef(t)=0$ and $Ef^2(t)=1$.
\item[(C2)] $f$ is almost surely at least three times continuously
differentiable with respect to $t$.
\item[(C3)] $T$ is a $d$-dimension Jordan measurable compact
subset of $R^d$.
\item[(C4)] The Hessian matrix of $C(t)$ at the origin is $-I$, where
$I$ is a $d\times d$ identity matrix.
\end{longlist}
Condition (C1) imposes unit variance. We will later study $\int_T
e^{\sigma f(t)}\,dt$ and treat $\sigma$ as an extra parameter. Condition
(C2) implies that $C(t)$ is at least 6 times
differentiable. In addition, the first, third and fifth
derivatives of~$C(t)$ evaluated at the origin are zero.
For any $\tilde f(t)$ such that $\Delta\tilde C(0)= \Sigma$ and
$|\Sigma|>0$, (C4) can always be achieved by an affine transformation on
the domain $T$ by letting $\tilde f(t) = f(\Sigma^{1/2}t)$ and
\[
\int_T e^{\sigma\tilde f(t)}\,dt=\int_T e^{\sigma f(\Sigma
^{1/2}t)}\,dt= |\Sigma|^{-1/2}\int_{\{s\dvtx\Sigma^{-1/2}s\in T\}}
e^{\sigma f(s)}\,ds,
\]
where for a symmetric matrix $\Sigma$ we let
$\Sigma^{1/2}$ be a symmetric matrix such that $\Sigma^{1/2}\Sigma
^{1/2}=\Sigma$.

For $\sigma>0$, let
%
%
\begin{equation}\label{Int}I_\sigma(A) = \int_A
e^{\sigma f(t)}\,dt
\end{equation}
for the Jordan measurable set $A\subset T$.
Of interest is
\[
P\bigl(I_\sigma(T)>b\bigr)
\]
as
$b\rightarrow\infty$. Equivalently, we may consider that the
variance of $f$ is $\sigma^2$. However, it is notionally
simpler to focus on a unit variance field and treat $\sigma$ as
a scale parameter.

We adopt the following notation. Let ``$\partial$'' and
``$\Delta$'' denote the gradient and Hessian matrix with
respect to $t$, and ``$\partial^2$'' denote the vector of second
derivatives with respect to $t$. The difference between ``$\Delta$''
and ``$\partial^2$'' is that, for a specific $t$, $\Delta f(t)$
is a $d\times d$ symmetric matrix whose upper triangle entries
are the elements of $\partial^2 f(t)$ which is a
$(d(d+1)/2)$-dimensional vector. Let $\partial_j$ denote the partial
derivative with respect to the $j$th component of $t =
(t_1,\ldots,t_d)$. We use similar notation for higher order
derivatives. For $b$ large enough, let $u$ be the unique
solution to
\[
(2\pi/\sigma)^{d/2}u^{-d/2}e^{\sigma u} =b.
\]
The uniqueness of $u$ is immediate by
noting that the left-hand side is monotone increasing with $u$
for all $u>d/(2\sigma)$. In addition, we use the following
notation and changes of variables:
\begin{eqnarray*}
\mu_{1}(t) &=&-(\partial_{1}C(t),\ldots,\partial_{d}C(t)),\\
\mu_{2}(t) &=& \bigl(\partial^2_{ii}C(t),i=1,\ldots,d; \partial^2
_{ij}C(t),i=1,\ldots,d,j=i+1,\ldots,d\bigr),\\
\mu_{02}^\top&=&\mu_{20}
=\mu_{2}(0),\qquad f(0) =u-w,\qquad \partial f(0) =y,\\
\partial^{2}f(0) &=& u\mu_{02}+z,\qquad
\Delta f(0) =-uI+\zz.
\end{eqnarray*}
The vector $\mu_{20}$ contains the
spectral moments of order two. Similar to $\Delta f(0)$ and
$\partial^2 f(0)$, $\zz$ is a symmetric matrix whose entries
consist of elements in $z$. We create different notation
because we will treat the second derivative of~$f$ as a matrix
when doing Taylor's expansion and as a vector when doing
integration. As stated in condition (C4), we have $\Delta C(0)
= -I$. Equivalently,~$\partial f(0)$ is a vector of independent
unit variance Gaussian r.v.'s. We plan to show
that in order to have $\int_T e^{f(t)}\,dt>b$, $\sup_T f(t)$
needs to reach a~level around $u$. The distance between
$f(0)$ and $u$ is denoted by $w$. In addition, since $(f(0),\partial^2
f(0))$ is
jointly independent of $\partial f(0)$, the distribution
of~$\partial f(0)$ is unaffected even if $f(0)$ reaches a high
level. Further, the covariance between~$f(0)$ and~$\partial^2
f(0)$ is $\mu_{20}$. Given $f(0)=u$, the conditional
expectation of $\partial^2 f(0)$ is $u\mu_{02}$. The distance between
$\partial^2 f(0)$ and this conditional expectation is denoted by
vector $z$.

A well-known result (see, e.g., Chapter 5.5 in
\cite{AdlTay07}) is that the joint distribution of
$(f(0),\partial^{2}f(0),\partial f(0),f(t))$ is multivariate
normal with mean zero and variance
\[
\pmatrix{
1 & \mu_{20} & 0 & C(t) \vspace*{2pt}\cr
\mu_{02} & \mu_{22} & 0 & \mu_{2}^{\top}(t) \vspace*{2pt}\cr
0 & 0 & I & \mu_{1}^{\top}(t) \vspace*{2pt}\cr
C(t) & \mu_{2}(t) & \mu_{1}(t) & 1},
\]
where $\mu_1(t)$, $\mu_2(t)$ and $\mu_{20}= \mu_{02}^\top$ is
defined previously.
The matrix $\mu_{22}$ is a $d(d+1)/2$ by $d(d+1)/2$ positive definite
matrix and contains the 4th-order
spectral moments arranged in an appropriate order. Conditional
on $f(0)= u-w$, $\partial f (0) = y$ and $\Delta f (0)
=-uI+\zz$, $f(t)$ is a continuous Gaussian random field with
conditional expectation
%
%
\begin{equation}\label{Et}
E(t)=(u-w,u\mu_{20}+z^{\top},y^{\top})
\pmatrix{
\Gamma^{-1} & 0 \cr
0 & I}
\pmatrix{C(t) \cr
\mu_{2}(t) \cr
\mu_{1}(t)},
\end{equation}
where
%
%
\begin{equation}\label{gamma}
\Gamma=\pmatrix{1 & \mu_{20} \cr
\mu_{02} & \mu_{22}}.
\end{equation}
Note that\vspace*{1pt} $u\mu_{20} +z$ is the vector
version of $-u I + \zz$. Therefore, conditional on
$(f(0),\partial f(0)^\top,\partial^2 f(0)^\top) = (u-w, y^\top,
u\mu_{20} + z^\top)$, we have representation
\[
f(t)=E(t) + g(t),\vadjust{\goodbreak}
\]
where $g(t)$ is a Gaussian random field
with mean zero and
\[
E(t) = E\bigl(f(t)| f(0) = u-w, \partial
f(0) = y, \partial^2 f(0) = u \mu_{02}+ z\bigr),
\]
whose form is given in (\ref{Et}). Since
$C(t)$ is six times differentiable, $E(t)$ is at least four
times differentiable. Using the form of $E(t)$ in (\ref{Et}) and
$\Gamma$ in (\ref{gamma}), after some tedious calculations, we have
%
%
\begin{eqnarray} \label{Derivative}
E(0) &=& u-w,\qquad \partial E(0)=y,\nonumber\\[-2pt]
\Delta E(0)&=&-uI+\zz,\qquad
\partial_{ijk}^{3}E(0)=y^{\top}\,\partial_{ijk}\mu_{1}(0),\\[-2pt]
\partial_{ijkl}^{4}E(0) &=& (u-w,u\mu_{20}+z^{\top})\Gamma
^{-1}\pmatrix{\partial_{ijkl}C(0) \cr
\partial_{ijkl}\mu_{2}(0)}.\nonumber
\end{eqnarray}
In order to obtain the above identities, we need the following facts.
The first, third and fifth derivatives of $C(t)$ evaluated at $0$ are
all zero. The first and second derivatives of $C(t)$ are contained in
$\mu_1(t)$ and $\mu_2(t)$. We also need to use the fact that
\[
\Gamma
^{-1}=\pmatrix{
\displaystyle \frac{1}{1-\mu_{20}\mu_{22}^{-1}\mu_{02}} & \displaystyle -\frac{\mu
_{20}\mu_{22}^{-1}
}{1-\mu_{20}\mu_{22}^{-1}\mu_{02}} \vspace*{3pt}\cr
\displaystyle -\frac{\mu_{22}^{-1}\mu_{02}}{1-\mu_{20}\mu
_{22}^{-1}\mu_{02}} & \displaystyle \mu_{22}^{-1}+\frac{\mu_{22}^{-1}\mu
_{02}\mu_{20}\mu_{22}^{-1}}{1-\mu
_{20}\mu_{22}^{-1}\mu_{02}}}.
\]
With the derivatives of $E(t)$, we can write
%
%
\begin{equation}\label{E}
E(t)=u-w+y^{\top}t+\tfrac{1}{2}t^{\top
}(-uI+\zz)t+g_{3}(t)+g_{4}(t) + R(t).
\end{equation}
If we let $t=(t_1,\ldots,t_d)$, then
%
%
\begin{equation}\label{G34}
g_3(t) = \frac1 6 \sum_{i,j,k}\partial^3_{ijk}E(0) t_i t_j t_k,\qquad
g_4(t) = \frac1 {24} \sum_{i,j,k,l}\partial^4_{ijkl}E(0) t_i t_j t_k
t_l,
\end{equation}
and $R(t)$ is the remainder term of the Taylor
expansion. The Taylor expansion of $E(t)$ is the same as
$f(t)$ for the first two terms because $g(t)$ is of order
$O(|t|^3)$. It is not hard to check that $\operatorname{Var}(g(t))\leq c|t|^6$
for some $c>0$ and $|t|$ small enough.

\subsection{Some related existing results}

For the comparison with the high excursion of $\sup_T f(t)$, we cite
one result for homogeneous random fields, which has been proved in more
general settings in many different ways. See, for instance,
\cite{Pit95,AW05,AdlTay07}. This result is also useful for the
proof of Theorem \ref{ThmM}. For comparison purpose, we only
present the result for the random fields
discussed in this paper.
\begin{proposition}\label{ThmPiterbarg}
Suppose Gaussian random field $f$ satisfies condi-\break tion~\textup{(C1)--(C4)}.
There exists a constant $G$ such that
\[
P\Bigl(\sup_{t\in T}f(t)>u\Bigr)=\bigl(1+o( 1)
\bigr)G\mes(T)u^{d}P\bigl(f(0)>u\bigr)
\]
as $u\rightarrow\infty$.\vadjust{\goodbreak}
\end{proposition}

We also present one existing result on the tail probability
approximation of the sum of correlated log-normal random
variables which provides intuitions on the analysis of $\int_T
e^{f(t)}\,dt$.
\begin{proposition}\label{PropLogN}
Let $X= (X_1,\ldots,X_n)$ be a multivariate Gaussian random variable
with mean $\mu$ and covariance matrix $\Sigma$, with
$\det(\Sigma)>0$. Then,
%
%
\begin{equation}\label{LogN}
P\Biggl(\sum_{i=1}^n e^{X_i} > b\Biggr) = \bigl(1+o(1)\bigr)
\sum_{i=1}^n P(e^{X_i} > b)
\end{equation}
as $b \rightarrow
\infty$.
\end{proposition}

The proof of this proposition can be found in \cite{AR08,FA09}.
This result implies that the large value of $\sum_{i=1}^n e^{X_i}$
is largely caused by one of the $X_i$'s being large. In the case that
$X_i$'s are independent, Proposition \ref{PropLogN} is a simple
corollary of the subexponentiality of log-normal distribution. Though the
$X_i$'s are correlated, asymptotically they are
tail-independent. The result presented in the next section can
be viewed as a natural generalization of Proposition
\ref{PropLogN}. Nevertheless, the techniques are quite
different from the following aspects. First, Proposition
\ref{PropLogN} requires $\Sigma$ to be nondegenerated.
For the continuous random fields, this is usually not true. As
shown in the analysis, we indeed need to study the sum of
random variables whose correlation converges to~$1$ when $b$
tends to infinity. Second, the approximation in (\ref{LogN}) is
for a sum of a~fixed number of random variables. The analysis
of the continuous field usually needs to handle the situation
that the number of random variables in a~sum grows to infinity
as $b\rightarrow\infty$. Last but not least, to obtain
approximations for $P(\int_T e^{f(t)}\,dt
>b)$, one usually needs to first obtain approximations for
$P(\int_{\Xi_\varepsilon} e^{f(t)}\,dt
>b)$ for some small domain $\Xi_\varepsilon\subset T$. We will
address all these issues in later sections.

For notation convenience, we write $a_u = O(b_u)$ if there exists a
constant $c>0$ independent of everything such that $a_u \leq
cb_u$ for all $u>1$, and $a_u = o(b_u)$ if $a_u/b_u \rightarrow
0$ as $u\rightarrow\infty$ and the convergence is uniform in
other quantities. We write $a_u = \Theta(b_u)$ if $ a_u
=O(b_u)$ and $b_u = O(a_u)$. In addition, we write $X_u =
o_p(1)$ if $X_u \stackrel{p}{\rightarrow}0$ as $u\rightarrow
\infty$.

\section{Main result} \label{SecMain}

The main theorem of this paper is stated as follows.
\begin{theorem} \label{ThmM}
Let $f$ be a Gaussian random field living on $T\subset R^d$
satisfying \textup{(C1)--(C4)}. Given $\sigma>0$, for $b$ large enough, $u$
is the unique solution to equation
%
%
\begin{equation}\label{u}
\biggl(\frac{2\pi}{\sigma}
\biggr)^{d/2}u^{-d/2}e^{\sigma u}=b.
\end{equation}
Then,
\[
P\biggl(\int_T e^{\sigma f(t)}\,dt>b\biggr) = \bigl(1+o(1)\bigr) H \mes(T) u^{d-1}
\exp(-u^2/2)
\]
as $b\rightarrow\infty$, where $\mes(T)$ is the Lebesgue
measure of $T$,
%
%
\begin{eqnarray}\label{H}
H&=&\frac{|\Gamma
|^{-1/2}\det(\mu_{22})^{1/2}e^{({(1/8) \mathbf1^\top
\mu_{22} \mathbf1+ (1/8) \sum_i
\partial^4_{iiii} C(0)})/{\sigma^2}}}{(2\pi
)^{(d+1)(d+2)/4}}\nonumber\hspace*{-35pt}\\[-8pt]\\[-8pt]
&&{} \times \int_{R^{d(d+1)/2}}
\exp\biggl\{-\frac1 2 \biggl[ B^\top B +\frac{
(\mu_{20}\mu_{22}^{-1/2}B + {
\mu_{20}\mathbf1}/({2\sigma}))^2}{1-\mu_{20}\mu
_{22}^{-1}\mu_{02}}
\biggr]\biggr\}\,dB ,\nonumber\hspace*{-35pt}
\end{eqnarray}
$\Gamma$ is defined in (\ref{gamma}), $\mu_{20}$, $\mu_{02}$, $\mu
_{22}$ are defined in the previous section and
\[
\mathbf{1}=(\underbrace{1,\ldots,1}_{d},\underbrace
{0,\ldots,0}_{d(d-1)/2})^{\top}.
\]
\end{theorem}
\begin{remark}
The integral in (\ref{H}) is clearly in an analytic form. We write it
as an integral because it arises naturally from the derivation.
\end{remark}
\begin{corollary}\label{CorSharp}
Let $f$ be a Gaussian random field living on $T\subset R^d$
satisfying \textup{(C1)--(C4)}. Adopting all the notation in Theorem \ref{ThmM},
let $\tilde b = b(2\pi/\sigma)^{-d/2}$ and
%
%
\begin{equation}\label{ut}
\tilde{u}=\frac{\log\tilde b}{\sigma}+\frac{d}{2\sigma}\log
\biggl( \frac{\log\tilde b}{\sigma}\biggr) +\biggl( \frac{d}{2}\biggr)
^{2}\frac{\log( ({\log\tilde b
})/{\sigma}) }{\sigma\log\tilde b}.
\end{equation}
Then,
\[
P\bigl(I_\sigma(T)>b\bigr) = \bigl(1+o(1)\bigr) H \mes(T) \tilde u^{d-1} \exp(-\tilde u^2/2).
\]
\end{corollary}
\begin{pf}
The result is immediate by the Taylor expansion on the left-hand side of equation
(\ref{u}) and note that $u-\tilde u =o(u^{-1})$.
\end{pf}

\begin{figure}[b]

\includegraphics{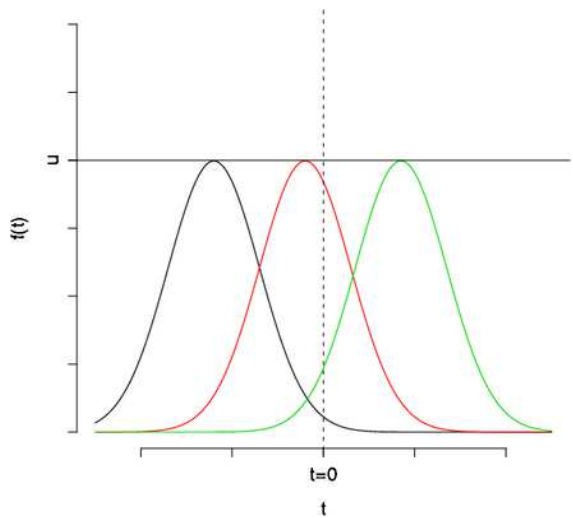}

\caption{One-dimensional example.} \label{FigShift}
\end{figure}

As we see, the asymptotic tail decaying rates of $\sup_T
f(t)$ and $\int_T e^{\sigma f(t)}\,dt$ take a very similar form.
More precisely,
\[
P\biggl(\int_T e^{\sigma f(t)}\,dt > b\biggr) =\Theta(1) P\Bigl(\sup_T
f(t) > u\Bigr)
\]
with $u$ and $b$ connected via (\ref{u}).
This fact suggests the following intuition on the tail
probability of $I_\sigma(T)$. First, the event
$\{I_\sigma(T)>b\}$ has substantial overlap with event $\{\sup_T
f(t)>u\}$. It has been shown by many studies mentioned before
that given $u$ sufficiently large $\{\sup_T f(t)>u\}$ is mostly caused
by just a single
$f(t^*)$ being large for some $t^*\in T$. Put these two facts
together, $\{I_\sigma(T)>b\}$ is mostly caused by $\{f(t^*)>u\}$,
for some $t^*\in T$ not too close to the boundary of $T$.
Therefore, conditional on $\{I_\sigma(T)>b\}$, the distribution
of $f(t)$ is very similar to the distribution conditional on
$\{\sup_T f(t)>u\}$. Of course, these two conditional
distributions are not completely identical. The difference will
be discussed momentarily. Now we perform some informal
calculation to illustrate the shape of $f(t)$ given $f(t^*)=u$.
Thanks to homogeneity, it is sufficient to study $t^*=0$. Given
$f(0)=u$,
\[
E\bigl(f(t)| f(0)=u\bigr)= uC(t).
\]
Since $C(t)$ is 6 times
differentiable, $\partial C(0) =0$ and $\Delta C(0) = -I$, we obtain
$E(f(t)| f(0)=u)\approx u - u\cdot t^\top t/2$. For the exact
Slepian model of the random field given that $f$ achieves a
local maximum at $t^*$ of level $u$, see \cite{ATW09}. Note
that for $b$ large,
\[
\int_T e^{\sigma u-(1/2)\sigma ut^\top t}\,dt > b
\]
is approximately equivalent to
\[
(2\pi/\sigma)^{d/2}u^{-d/2}e^{\sigma u} >b.
\]
In Theorem \ref{ThmM}, this is
exactly how $u$ is defined. As shown in Figure \ref{FigShift},
the three curves are $\exp\{E(f(t)|f(t^*)=u)\}$ for different
$t^*$'s. Given that \mbox{$\{\sup_T f (t)>u\}$}, these three curves are equally
likely to occur.

%
%

Second, as mentioned before, the conditional distributions of
$f(t)$ are different given $\{I_\sigma(T) > b\}$ or $\{\sup_T
f(t)>u\}$. This is why the two constants in Theorem~\ref{ThmM} ($H$)
and Proposition \ref{ThmPiterbarg} ($G$) are different. The
difference is due to the fact that the symmetric difference
between $\{\sup_T f(t)>u\}$ and $\{\int_T e^{f(t)}\,dt >b\}$ is
substantial though their overlap is significant too.
Consider the following situation that contributes to the
difference. $\sup_T f(t)$ is slightly less than $u$ [e.g.,
by a magnitude of $O(u^{-1})$]. For this case,
$I_\sigma(T)$ still has a large chance to be greater than $b$.
For this sake we will need to consider the contribution of $\Delta
f(0)$. As is shown in the technical proof, if $t^*=\arg\sup f(t) =0$
and $\partial
f(0)=0$, then a sufficient and necessary condition for
$I_\sigma(T)>b$ is that
\[
f(0) + \frac1 {2\sigma}\operatorname{Tr} \bigl(
u^{-1}\Delta f(0)+I \bigr)
>u + o(u^{-1}),
\]
where $\mathrm{Tr}$ denotes the trace of a squared
matrix. Note that conditional on $f(0)=u$, $E(\Delta f(0) |f(0)
= u) = -uI$. Therefore, $\Delta f(0)+uI$ is of size $O(1)$. One
well-known result is that the trace of a symmetric matrix is
the sum of its eigenvalues. Let $\lambda_i$ be the eigenvalues
of $u^{-1}\Delta f(0) +I = \zz/u$. Then, the sufficient and necessary
condition is
translated to $f(0)+\frac1 {2\sigma}\sum_{i=1}^d \lambda_i
>u$. This also suggests that, conditional on $I_\sigma(T) > b$,
$w= f(0) - u$ is of size~$O(u^{-1})$. This forms the
intuition behind the proof of Theorem \ref{PropSmall}.

The proof of Theorem \ref{ThmM} consists of two steps presented in
Sections \ref{SecStep1} and \ref{SecStep2}, respectively. Each
of the two steps is summarized as one theorem.

\subsection{Step 1}\label{SecStep1}

Construct a cover of $T$, $\{A_1,\ldots,A_n\}$, such that
$T\subset\bigcup_{i=1}^n A_i$. Each $A_i$ is a closed square,
$\mes(A_i\cap A_j) =0$ for $i\neq j$. Because $T$ is Jordan
measurable, as $\sup_i \mes(A_i)\rightarrow0$,
$\mes(\bigcup_{i=1}^n A_i)-\mes(T) \rightarrow0$. To simplify the
analysis, we make each $A_i$ of identical shape and let $A_i =
\{t_i + s \dvtx\allowbreak s\in[0,a_b]^d\}$. The size of the partition $n$
and choice of $a_b$ depend on the threshold $b$. The first step
analysis involves computing the integral $p_i \triangleq
P(\int_{A_i} e^{\sigma f(t)}\,dt>b)$. Because $f$ is homogeneous,
it is sufficient to study~$p_1$.

The basic strategy to approximate $p_1$ is as follows. Because
$f$ is at least three times differentiable, the first and
second derivatives are almost surely well defined. Without loss
of generality, we assume that $0\in A_1$. Conditional on
$(f(0),\partial f(0),\Delta f(0))$, $f(t)= f(0) + \partial
f(0)^\top t + \frac1 2 t^\top\Delta f(0) t + g_3(t) +
g_4(t)+R(t)+ g(t)$, where $g(t)$ is a Gaussian field with mean
zero and variance of order $O(|t|^6)$. Then,
%
%
\begin{eqnarray}\label{p1}\qquad
p_1 &=& P\biggl(\int_{A_1} e^{\sigma f(t)}\,dt >b\biggr)\nonumber\\
&=&\int h(w,y,z) \nonumber\\[-8pt]\\[-8pt]
&&\hphantom{\int }
{}\times P\biggl(\int_{A_1} e^{\sigma[u-w + y^\top t+ (1/2)
t^\top(-uI+\zz) t + g_3(t) + g_4(t) + R(t)+g(t)]}\,dt
\nonumber\\
&&\hphantom{\int {}\times P\biggl(}
{}>b \biggr)\,dw\,dy\,dz,\nonumber
\end{eqnarray}
where $h(w,y,z)$ is the density function of
$(f(0),\partial f(0),\Delta f(0))$ evaluated at
$(u-w,y,u\mu_{20}+z)$, which is a multivariate Gaussian random
vector. Let~$u$ be defined in (\ref{u}). There exists a $\delta
>0$ (small enough) such that if we let~$A_i$'s be squares of size
$\varepsilon= O(1)u^{-1/2+\delta}$ and, hence,
$n=O(1)\mes(T)u^{d/2-d\delta}$, the asymptotics of $p_i$ can be
derived by repeatedly using Taylor's expansion and evaluating
the integral on the right-hand side of (\ref{p1}). The main
result of this step is presented as follows. It establishes a
similar result to that of Theorem \ref{ThmM} but within a much smaller
domain.
\begin{theorem}\label{PropSmall}
Let $f$ be a Gaussian random field living in $T$ satisfying
conditions \textup{(C1)--(C4)}. Let $A_1 = \Xi_\varepsilon= \{t \dvtx
|t|_\infty< \varepsilon\}$, where $|t|_\infty= {\max_i} |t_i|$.
Let $u$ and $H$ be defined in Theorem \ref{ThmM}. Without loss
of generality, assume $\Xi_\varepsilon\subset T$ with
$\varepsilon= \kappa u^{\delta- 1/2}$ for some $\delta$ small
enough and $u$ large enough. Then, for any $\kappa>0$
\[
p_1=p(\Xi_\varepsilon)=P\biggl(\int_{\Xi_{\varepsilon}}e^{\sigma f(t)}\,dt
> b\biggr) = \bigl(1+o(1)\bigr)H
\mes(\Xi_\varepsilon)u^{d-1}e^{-u^2/2}
\]
as $b \rightarrow
\infty$.
\end{theorem}

The proof of this theorem is in Section \ref{sec_proof1}. We
will then choose each $A_i$ to be of the same shape as
$\Xi_\varepsilon$. Then, all the $p_i$'s are identical.

\subsection{Step 2}\label{SecStep2}

The second step is to show that with the particular choice
of~$A_i$ in the first step, $P(\int_{T} e^{\sigma f(t)}\,dt>b) =
(1+o(1))\sum_{i=1}^n p_i$. We first present the main result of
the second step.
\begin{theorem}\label{PropInd} Let $f$ be a Gaussian random
field satisfying conditions~\textup{(C1)--(C4)} and $\varepsilon$ be chosen in
Theorem \ref{PropSmall}. Let $\kk\in\mathbb Z^d$ and
$\Xi_{\varepsilon,\kk}= 2\kk\varepsilon+
\Xi_\varepsilon$. Further, let $\mathcal C^-
=\{\kk\dvtx\Xi_{\varepsilon,\kk} \subset T\}$ and $\mathcal C^+
=\{\kk\dvtx\Xi_{\varepsilon,\kk} \cap T\neq\varnothing\}$, then
\[
P\biggl(I_\sigma\biggl(\bigcup_{\kk\in\mathcal C^+} \Xi_{\varepsilon,\kk}\biggr)>b\biggr)
= \bigl(1+o(1)\bigr)\sum_{\kk\in\mathcal C^+}
P\bigl(I_\sigma(\Xi_{\varepsilon,\kk})>b\bigr)
\]
and
\[
P\biggl(I_\sigma\biggl(\bigcup_{\kk\in\mathcal C^-} \Xi_{\varepsilon,\kk}\biggr)>b\biggr)
= \bigl(1+o(1)\bigr)\sum_{\kk\in\mathcal C^-}
P\bigl(I_\sigma(\Xi_{\varepsilon,\kk})>b\bigr).
\]
\end{theorem}

We consider
\[
I_\sigma\biggl(\bigcup_{\kk\in\mathcal C^+}
\Xi_{\varepsilon,\kk}\biggr)= \sum_{\kk\in\mathcal C^+}
I_\sigma(\Xi_{\varepsilon,\kk})
\]
as a sum of finitely many
dependently and identically distributed random variables. The
conclusion of the above theorem implies that the tail
distribution of\vadjust{\goodbreak} the sum of these dependent variables
exhibits the so-called ``one big jump'' feature---the high
excursion of the sum is mainly caused by just one component
being large. This result is similar to that of the sum of
correlated log-normal r.v.'s. Nevertheless, the gap between the
analyses of finite sum and integral is substantial because the
correlation between fields in adjacent squares tends to $1$. For
finite sums, the correlation is always bounded away from $1$.
The key step in the proof of Theorem \ref{PropInd} is that
the $\varepsilon$ defined in Theorem \ref{PropSmall},
though tends to zero as $b\rightarrow\infty$, is large enough
such that the one-big-jump principle still applies. We will connect the
event of high excursion of
$I_\sigma(\bigcup_{\kk\in\mathcal C^+}
\Xi_{\varepsilon,\kk})$ to the high excursion of
$\sup_{\bigcup_{\kk\in\mathcal C^+} \Xi_{\varepsilon,\kk}}f
(t)$ and apply existing results on the bound on the supremum of
Gaussian random fields. A short list of recent related
literature on the ``one-big-jump'' principle and multivariate
Gaussian random variables is \cite{FA09,DDS08,AR08}.

With the preparation of the two steps, we are ready to
present the proof for Theorem \ref{ThmM}.
\begin{pf*}{Proof of Theorem \ref{ThmM}}
From Theorem \ref{PropSmall},
\begin{eqnarray*}\sum_{\kk
\in\mathcal C^+} P\bigl(I_\sigma(\Xi_{\varepsilon,\kk})>b\bigr) &=&
\bigl(1+o(1)\bigr) H \mes\biggl(\bigcup_{\kk\in
\mathcal{C}^{+}}\Xi_{\varepsilon,\kk}\biggr) u^{d-1}e^{-u^2/2},\\
\sum_{\kk\in\mathcal C^-}
P\bigl(I_\sigma(\Xi_{\varepsilon,\kk})>b\bigr) &=& \bigl(1+o(1)\bigr) H
\mes\biggl(\bigcup_{\kk\in\mathcal{C}^{-}}\Xi_{\varepsilon
,\kk}\biggr) u^{d-1}e^{-u^2/2}.
\end{eqnarray*}
Therefore,
thanks to Theorem \ref{PropInd},
\begin{eqnarray*}
P\bigl(I_\sigma(T)>b\bigr)&\geq& P\biggl(I_\sigma\biggl(\bigcup_{\kk\in\mathcal{C}^{-}}\Xi
_{\varepsilon,\kk}\biggr)>b\biggr)\\
&\geq&\bigl(1+o(1)\bigr) H \mes\biggl(\bigcup_{\kk\in
\mathcal{C}^{-}}\Xi_{\varepsilon,\kk}\biggr)
u^{d-1}e^{-u^2/2};
\end{eqnarray*}
similarly,
\begin{eqnarray*}
P\bigl(I_\sigma(T)>b\bigr)&\leq& P\biggl(I_\sigma\biggl(\bigcup_{\kk\in\mathcal{C}^{+}}\Xi
_{\varepsilon,\kk}\biggr)>b\biggr)\\
&\leq&\bigl(1+o(1)\bigr) H \mes\biggl(\bigcup_{\kk\in
\mathcal{C}^{+}}\Xi_{\varepsilon,\kk}\biggr)
u^{d-1}e^{-u^2/2}.
\end{eqnarray*}
Jordan measurability of $T$ implies that
\[
\mes\biggl(\bigcup_{\kk\in
\mathcal{C}^{+}}\Xi_{\varepsilon,\kk}\biggr) - \mes\biggl(\bigcup_{\kk\in
\mathcal{C}^{-}}\Xi_{\varepsilon,\kk}\biggr) \rightarrow0{+}.
\]
Therefore,
\[
P\bigl(I_\sigma(T)>b\bigr) = \bigl(1+o(1)\bigr) H \mes(T)u^{d-1}e^{-u^2/2}.
\]
\upqed\end{pf*}

\section{\texorpdfstring{Proof for Theorem \protect\ref{PropSmall}}{Proof for Theorem 2}}
\label{sec_proof1}

In this section we present the proof of Theorem~\ref{PropSmall}. We arrange all the lemmas
and their proofs in the \hyperref[app]{Appendix}.
\begin{pf*}{Proof of Theorem \ref{PropSmall}}
We evaluate the probability by conditioning on $(f(0),\partial
f(0), \partial^2 f(0))$,
%
%
\begin{eqnarray} \label{IntCond}
p(\Xi_{\varepsilon}) &=&P\biggl( \int_{\Xi_{\varepsilon
}}e^{\sigma f(t)}\,dt>b\biggr) \nonumber\\
&=&\int_{\mathcal R} h(w,y,z)\nonumber\\
&&\hspace*{11pt}{}\times P\biggl( \int_{\Xi
_{\varepsilon}}e^{\sigma f(t)}\,dt>b\Big\vert
f(0)=u-w,\partial f(0)=y,\\
&&\hspace*{146.6pt}\partial^{2}f(0)=u\mu
_{02}+z\biggr) \,dw\,dy\,dz \nonumber\\
&=&\int_{\mathcal R} h(w,y,z)P\biggl( \int_{\Xi
_{\varepsilon}}e^{\sigma E(t)+\sigma g(t)}\,dt>b\biggr) \,dw\,dy\,dz,\nonumber
\end{eqnarray}
where $\mathcal R = R^{(d+1)(d+2)/2}$ and $h(w,y,z)$ is the
density function of $(f(0), \partial f(0)$, $\partial^2 f(0))$
evaluated at $( u-w, y, u\mu_{02}+z)$. Now we take a closer
look at the integrand inside the above integral.\vspace*{2pt} Conditional on
$f(0)=u-w,\partial f(0)=y$, $\partial^{2}f(0)=u\mu_{02}+z$,
\begin{eqnarray*}
I_\sigma(\Xi_\varepsilon)
&=&\int_{|t|_{\infty
}<\varepsilon}e^{\sigma E(t)+\sigma g(t)}\,dt \\
&=&\int_{|t|_{\infty}<\varepsilon}\exp\biggl\{\sigma
\biggl[u-w+y^{\top}t+\frac{1}{2}t^{\top
}(-uI+\zz)t\\
&&\hphantom{\int_{|t|_{\infty}<\varepsilon}\exp\biggl\{}
\hspace*{17pt}{}+g_{3}(t)+g_{4}(t)+ R(t)+g(t)\biggr]\biggr\}\,dt \\
&=&\det(uI-\zz)^{-1/2}\\
&&{}\times\int_{|(uI-\zz)^{-1/2}t|_{\infty
}<\varepsilon}e^{\sigma\{u-w+
({1/2})y^{\top}(uI-\zz)^{-1}y\}} \\
&&{}\times\exp\biggl\{\sigma\biggl[-\frac{1}{2}\bigl(t-(uI-\zz)^{-1/2}y\bigr)^{\top}
\bigl(t-(uI-\zz)^{-1/2}y\bigr)\\
&&\hspace*{45.2pt}
{} +g_{3}\bigl((uI-\zz)^{-1/2}t\bigr)
+g_{4}\bigl((uI-\zz)^{-1/2}t\bigr)\\
&&\hspace*{69.6pt}
{}+R\bigl((uI-\zz)^{-1/2}t\bigr)+g\bigl((uI-\zz
)^{-1/2}t\bigr)
\biggr]\biggr\}\,dt.
\end{eqnarray*}
For the\vspace*{1pt} second equality, we plugged in (\ref{E}). For the
last step, we first change the variable from $t$ to
$(uI-\zz)^{1/2}t$ and then write the exponent in a quadratic
form of~$t$. We write the term inside the exponent without the
factor $\sigma$ as
%
%
\begin{eqnarray}\label{J}\qquad
J(t)&=&-\tfrac{1}{2}\bigl(t-(uI-\zz)^{-1/2}y\bigr)^{\top
}\bigl(t-(uI-\zz)^{-1/2}y\bigr)+g_{3}\bigl((uI-\zz)^{-1/2}t\bigr)\nonumber\\[-8pt]\\[-8pt]
&&{}+g_{4}\bigl((uI-\zz
)^{-1/2}t\bigr)+R\bigl((uI-\zz)^{-1/2}t\bigr),\nonumber
\end{eqnarray}
which is asymptotically a quadratic form. But, as is shown
later, $g_3$ and $g_4$ terms do play a role in the calculation.
Also, it is useful to keep in mind that $J(t)$ depends on $y$ and $z$.
Hence, we can write
\begin{eqnarray*}
I_\sigma(\Xi_\varepsilon) &=&
\int_{|t|_\infty<\varepsilon}e^{\sigma f(t)}\,dt\\
&=&\det
(uI-\zz)^{-1/2}e^{\sigma\{u-w+({1/2})y^{\top
}(uI-\zz)^{-1}y\}}\\
&&{}\times\int_{|(uI-\zz)^{-1/2}t|_\infty<\varepsilon
}e^{\sigma J(t)+\sigma g((uI-\zz)^{-1/2}t)}\,dt.
\end{eqnarray*}
Let
%
%
\begin{equation}\label{H0}
e^{H_{0}}=\int_{R^d} e^{-({\sigma/2})t^{\top}t}\,dt
= (2\pi/\sigma)^{d/2}.
\end{equation}
Let $u$ solve
%
%
\begin{equation}\label{u0}
u^{-d/2}e^{\sigma u+H_{0}}=b.
\end{equation}
Then,
\[
\int_{\Xi_\varepsilon}e^{\sigma f(t)}\,dt>b,
\]
if and only if
\begin{eqnarray*}
\hspace*{-4pt}&&\det(uI-{\mathbf z})^{-1/2}e^{\sigma
\{u-w+({1/2})y^{\top}(uI-{\mathbf z})^{-1}y\}}\int_{|(uI-{\mathbf
z})^{-1/2}t|_{\infty
}<\varepsilon}e^{\sigma J(t)+\sigma g((uI-{\mathbf z})^{-1/2}t)}\,dt\\
\hspace*{-4pt}&&\qquad>u^{-d/2}e^{\sigma u+H_{0}}.
\end{eqnarray*}
We take the logarithm on both sides and rewrite the above inequality
and have
%
%
\begin{eqnarray}
0 &<&\frac{\sigma}{2}y^{\top}(uI-\zz)^{-1}y-\sigma
w-\frac{1}{2}\log\det(I-u^{-1}\zz) \\
&&{}+\log
\int_{|(uI-\zz)^{-1/2}t|_\infty<\varepsilon
}e^{\sigma J(t)}\,dt-H_{0}\nonumber\\
&&{}+\log E\exp\bigl(\sigma g\bigl((uI-\zz
)^{-1/2}S\bigr)\bigr)\nonumber\\
\label{Remainder}
&=& A(w,y,z) + \log E\exp\bigl(\sigma
g\bigl((uI-\zz)^{-1/2}S\bigr)\bigr),
\end{eqnarray}
where $S$ is a random variable on the region that $|(uI-\zz
)^{-1/2}S|_\infty\leq\varepsilon$ with density proportional to
$e^{\sigma
J(s)}$ and
%
%
\begin{eqnarray}\label{DefA}
A(w,y,z)&=&\frac{\sigma}{2}y^{\top}(uI-\zz)^{-1}y-\sigma
w-\frac{1}{2}\log\det(I-u^{-1}\zz)\nonumber\\[-9pt]\\[-9pt]
&&{} + \log
\int_{|(uI-\zz)^{-1/2}t|<\varepsilon}e^{\sigma
J(t)}\,dt-H_{0}.\nonumber
\end{eqnarray}

Thanks to Lemma \ref{LemEx1}, we only need to consider the set
that
%
%
\begin{eqnarray}\label{DefL}
\mathcal L&=&\{\vert f( 0) -u\vert\leq
u^{2\delta
+\varepsilon_{0}},\nonumber\\[-9pt]\\[-9pt]
&&\hphantom{\{}\vert\partial f( 0) \vert\leq u^{{1
}/{2}+\delta+\varepsilon_{0}},\vert\partial
^{2}f(0)-u\mu_{20}\vert\leq u^{{1/2}+\varepsilon
_{0}}\}.\nonumber
\end{eqnarray}
Also, by abusing notation, we write
%
%
\begin{equation}\label{L}\mathcal L=\{\vert w\vert\leq
u^{2\delta+\varepsilon_{0}},\vert y \vert\leq u^{{1
}/{2}+\delta+\varepsilon_{0}},\vert z\vert\leq
u^{{1}/{2}+\varepsilon_{0}}\}.
\end{equation}

Lemma \ref{LemDensity} gives the form of $h(w,y,z)$. We plug in the
results in Lemmas~\ref{LemEx1} and~\ref{LemDensity},
%
%
\begin{eqnarray}\label{ProbS}
p(\Xi_{\varepsilon})
&=&\int_{\mathcal{R}}h(w,y,z)\nonumber\\[-2pt]
&&\hphantom{\int_{\mathcal{R}}}{}\times P\bigl(A(w,y,z)+\log
E\exp\bigl(\sigma g\bigl((uI- \zz)^{-1/2}S\bigr)\bigr)>0\bigr)\,dw\,dy\,dz \nonumber\\[-2pt]
&=&o(1)u^{-\alpha
}e^{-u^{2}/2}\nonumber\\[-0.4pt]
&&{}+\int_{\mathcal{L}}h(w,y,z)\nonumber\\[-2pt]
&&\hphantom{{}+\int_{\mathcal{L}}}
{}\times P\bigl(A(w,y,z)+\log
E\exp\bigl( \sigma g\bigl((uI-\zz)^{-1/2}S\bigr)\bigr)>0 \bigr)\,dw\,dy\,dz \nonumber\\[-2pt]
&=&o(1)u^{-\alpha}e^{-u^{2}/2} \nonumber\\[-2pt]
&&{} +\frac{1}{(2\pi)^{(d+1)(d+2)/4}}|\Gamma
|^{-1/2}\\[-2pt]
&&\hspace*{11pt}{}\times\int_{\mathcal{L}}P\bigl(A(w,y,z) +\log
E\exp\bigl(\sigma g\bigl((uI-\zz)^{-1/2}S\bigr)\bigr)>0\bigr)\nonumber\\[-2pt]
&&\hspace*{34.2pt}{}\times\exp\biggl\{-\biggl[\frac{1}{2}u^{2}+\frac u \sigma
A(w,y,z)\nonumber\\[-2pt]
&&\hspace*{81.3pt}{}+\frac{1}{2}y^{\top
}\bigl(I-(I-{\mathbf z}/u)^{-1}\bigr)y \nonumber\\[-2pt]
&&\hspace*{81.3pt}{} +\frac{1}{2}\frac{(w+\mu_{20}\mu_{22}^{-1}z)^{2}}{1-\mu
_{20}\mu_{22}^{-1}\mu_{02}}+\frac{1}{2}z^{\top}\mu
_{22}^{-1}z\nonumber\\[-2pt]
&&\hspace*{81.3pt}{}+\frac{u}{2\sigma}\log
\det(I-u^{-1}\zz) \nonumber\\[-2pt]
&&\hspace*{81.3pt}{} -\frac u\sigma\log\int_{|(uI-\mathbf
z)^{-1/2}t|<\varepsilon}e^{\sigma J(t)}\,dt+\frac u\sigma
H_{0}\biggr]\biggr\} \,dw\,dy\,dz.\nonumber\vadjust{\goodbreak}
\end{eqnarray}
We define
%
%
\begin{eqnarray}\label{I}
\mathcal I &=&\frac{1}{2}u^{2}+\frac u\sigma
A(w,y,z)+\frac{1}{2}y^{\top}\bigl(I-(I-
\zz/u)^{-1}\bigr)y \nonumber\\
&&{} +\frac{1}{2}\frac{(w+\mu_{20}\mu_{22}^{-1}z)^{2}}{1-\mu
_{20}\mu_{22}^{-1}\mu_{02}}+\frac{1}{2}z^{\top}\mu
_{22}^{-1}z+\frac{u}{2\sigma}\log
\det(I-u^{-1}\zz) \\
&&{} -\frac u\sigma\log
\int_{|(uI-\mathbf{z})^{-1/2}t|<\varepsilon}e^{\sigma J(t)}\,dt+\frac u
\sigma H_{0},\nonumber
\end{eqnarray}
and proceed with some tedious algebra to write $\mathcal
I$ in a friendly form for integration. First notice that
\[
(I-u^{-1}\zz)^{-1} = \sum_{k=0}^\infty u^{-n}\zz^n.
\]
Plug this
into the third term of $\mathcal I$ and obtain
\begin{eqnarray*}
\mathcal{I} &=&\frac{1}{2}u^{2}+\frac u\sigma
A(w,y,z)-\bigl(1+O(|z|/u)\bigr)\frac{u^{-1}}{2}
y^{\top}\zz y \\
&&{} +\frac{1}{2}\frac{(w+\mu_{20}\mu_{22}^{-1}z)^{2}}{1-\mu
_{20}\mu_{22}^{-1}\mu_{02}}+\frac{1}{2}z^{\top}\mu
_{22}^{-1}z+\frac{u}{2\sigma}\log
\det(I-u^{-1}{\mathbf z}) \\
&&{} -\frac u\sigma\log
\int_{|(uI-\mathbf{z})^{-1/2}t|<\varepsilon}e^{\sigma J(t)}\,dt+\frac
u\sigma H_{0}.
\end{eqnarray*}

\subsection*{\texorpdfstring{Situation 1 of Lemma \protect\ref{LemJ}}
{Situation 1 of Lemma 5}} Adopt the
notation in Lemmas \ref{lemG} and \ref{LemJ}. Note that
according to the definition of $Y$ in Lemma
\ref{lemG} that
\begin{eqnarray*}
Y &=&(y_{i}^{2},i=1,\ldots,d,2y_{i}y_{j},1\leq i<j\leq d)^{\top}, \\
\mathbf{1}&=&(\underbrace{1,\ldots,1}_{d},
\underbrace{0,\ldots,0}_{d(d-1)/2})^{\top},
\end{eqnarray*}
we obtain that
\[
y^{\top}\zz y= Y^\top z.
\]
We plug\vspace*{1pt} in results of
Lemmas \ref{lemG} and \ref{LemJ}. First, considering the first
situation in Lemma \ref{LemJ}, that is, $\mathcal L_1=\mathcal L
\cap\{|(uI-\zz)^{-1/2}y|_\infty\leq\kappa
u^{\delta}-u^{\delta/2}\}$, we have
\begin{eqnarray*}
\mathcal{I} &=&\frac{1}{2}u^{2}+\frac u \sigma
A(w,y,z)-\bigl(1+O(|z|/u)\bigr)\frac{u^{-1}}{2}%
Y^{\top}z \\
&&{} +\frac{1}{2}\frac{(w+\mu_{20}\mu_{22}^{-1}z)^{2}}{1-\mu
_{20}\mu_{22}^{-1}\mu_{02}}+\frac{1}{2}z^{\top}\mu
_{22}^{-1}z+\frac{u}{2\sigma}\log
\det(I-u^{-1}\zz) \\
&&{} +\frac{1}{8}(u^{-1}Y+\mathbf{1}/\sigma)^{\top}\mu
_{22}(u^{-1}Y+\mathbf{1}/\sigma)-
\frac{1}{8\sigma^2}\mathbf{1}^{\top}\mu_{22}\mathbf{1} \\
&&{} +\frac u\sigma H_{0}-\frac u \sigma
H\bigl((uI-\mathbf{z})^{-1/2}y,(uI-\mathbf{z})^{1/2},\varepsilon\bigr)\\
&&{} -\frac{1}{8\sigma^2}\sum_{i}\partial_{iiii}C(0)+o(1).
\end{eqnarray*}
Also, it is useful to keep in mind that $\mathbf1$ is NOT a
vector of $1$'s. The next step is to plug in the result of
Lemma \ref{LemDet} and replace the $\log\det(I-u^{-1}\zz)$
term by
\begin{eqnarray*}
&&-u^{-1}\operatorname{Tr}(\zz) + \tfrac1 2 u^{-2} \mathcal E^2 (\zz) +o(u^{-1})\\[2pt]
&&\qquad=
-u^{-1} \mathbf1^\top z + \tfrac1 2 u^{-2} \mathcal E^2
(\zz) +o(u^{-1})
\end{eqnarray*}
and obtain
\begin{eqnarray*}
\mathcal{I} &=&\frac{1}{2}u^{2}+\frac u \sigma
A(w,y,z)-\bigl(1+O(|z|/u)\bigr)\frac{u^{-1}}{2}
Y^{\top}z \\[2pt]
&&{} +\frac{1}{2}\frac{(w+\mu_{20}\mu_{22}^{-1}z)^{2}}{1-\mu
_{20}\mu_{22}^{-1}\mu_{02}}+\frac{1}{2}z^{\top}\mu
_{22}^{-1}z-\frac{1}{2\sigma}
\mathbf{1}^{\top}z + \frac1 {4\sigma u} \mathcal E^2 (\zz) \\[2pt]
&&{} +\frac{1}{8}(u^{-1}Y+\mathbf{1}/\sigma)^{\top}\mu
_{22}(u^{-1}Y+\mathbf{1}/\sigma)-
\frac{1}{8\sigma^2}\mathbf{1}^{\top}\mu_{22}\mathbf{1} \\[2pt]
&&{} +\frac u \sigma H_{0}-\frac u \sigma
H\bigl((uI-\mathbf{z})^{-1/2}y,(uI-\mathbf{z})^{1/2},\varepsilon\bigr)\\[2pt]
&&{}-\frac{1}{8\sigma^2}\sum_{i}\partial_{iiii}C(0)+o(1).
\end{eqnarray*}
Then, we group the terms $-(1+O(|z|/u))\frac{u^{-1}}{2}
Y^{\top}z $ and $-\frac{1}{2\sigma}
\mathbf{1}^{\top}z$ and leave the $O(|z|/u)$ to the end and have
\begin{eqnarray*}
\mathcal I&=&\frac{1}{2}u^{2}+\frac u \sigma
A(w,y,z)+\frac{1}{2}\frac{(w+\mu_{20}\mu
_{22}^{-1}z)^{2}}{1-\mu_{20}\mu_{22}^{-1}\mu_{02}} \\[2pt]
&&{} +\frac{1}{2}z^{\top}\mu
_{22}^{-1}z-\frac{1}{2}(u^{-1}Y+\mathbf{1}/\sigma)^{\top
}z\\[2pt]
&&{}+\frac{1}{8}(u^{-1}Y+\mathbf{1}/\sigma)^{\top}\mu
_{22}(u^{-1}Y+\mathbf{1}/\sigma) \\[2pt]
&&{} +\frac u\sigma H_{0}-\frac u\sigma H\bigl((uI-\mathbf
{z})^{-1/2}y,(uI-\mathbf{z})^{1/2},\varepsilon\bigr) \\[2pt]
&&{} -\frac{1}{8\sigma^2}\mathbf{1}^{\top}\mu
_{22}\mathbf{1} -\frac{1}{8\sigma^2}
\sum_{i}\partial_{iiii}C(0)+o(1)\\[2pt]
&&{}+O(u^{-2}|z|^2|y|^2) + O(u^{-1}
\mathcal E^2 (\zz)).
\end{eqnarray*}
Note that second and third lines in the above display is in fact in
a~quadratic form. We then have
\begin{eqnarray*}
\mathcal I&=&\frac{1}{2}u^{2}+\frac u \sigma
A(w,y,z)+\frac{1}{2}\frac{(w+\mu_{20}\mu
_{22}^{-1}z)^{2}}{1-\mu_{20}\mu_{22}^{-1}\mu_{02}} \\[-2pt]
&&{} +\frac{1}{2}\biggl[\mu_{22}^{-1/2}z-\frac{1}{2}\mu
_{22}^{1/2}(u^{-1}Y+\mathbf{1}/\sigma)\biggr]^{\top}\\[-2pt]
&&\hspace*{11pt}{}\times\biggl[\mu
_{22}^{-1/2}z-\frac{1}{2}\mu
_{22}^{1/2}(u^{-1}Y+\mathbf{1}/\sigma)\biggr]
\\[-2pt]
&&{} +\frac u \sigma H_{0}-\frac u \sigma H\bigl((uI-\mathbf
{z})^{-1/2}y,(uI-\mathbf{z})^{1/2},\varepsilon\bigr)\\[-2pt]
&&{}-\frac{1}{8\sigma^2}\mathbf{1}^{\top}\mu
_{22}\mathbf{1}-\frac{1}{8\sigma^2} \sum_{i}\partial
_{iiii}C(0)+o(1)\\[-2pt]
&&{}+O(u^{-2}|z|^2|y|^2) + O(u^{-1} \mathcal E^2
(\zz)).
\end{eqnarray*}
Now, consider another change of variable,
%
%
\begin{eqnarray}\label{change}
A&=&A(w,y,z),\nonumber\\[-9pt]\\[-9pt]
B&=&\mu
_{22}^{-1/2}z-\tfrac{1}{2}\mu_{22}^{1/2}(u^{-1}Y+\mathbf{1}/\sigma
),\qquad
y=y.\nonumber
\end{eqnarray}
Then, by noting that $\mu_{20}$ is a row vector in which the first $d$
entries are $-1$'s and the rest are $0$'s, we have
\[
w+\mu_{20}\mu_{22}^{-1}z=-\frac{A}{\sigma}+\mu_{20}\mu
_{22}^{-{1/2%
}}B+\frac{1}{2\sigma}\mu_{20}\mathbf{1}+o(1).
\]
Therefore, we have
\begin{eqnarray*}
\mathcal{I} &=&\frac{1}{2}u^{2}+\frac{u}{\sigma}A+\frac
{1}{2}B^{\top}B\\[-2pt]
&&{}+
\frac{1}{2}\frac{(-{A}/{\sigma}+\mu_{20}\mu_{22}^{-
{1/2}}B+%
({1}/({2\sigma}))\mu_{20}\mathbf{1}+o(1))^{2}}{1-\mu_{20}\mu
_{22}^{-1}\mu_{02}}\\[-2pt]
&&{}+\frac{u}{\sigma}H_{0}
-\frac{u}{\sigma}H\bigl((uI-\mathbf
{z})^{-1/2}y,(uI-%
\mathbf{z})^{1/2},\varepsilon\bigr)\\[-2pt]
&&{} -\frac{1}{8\sigma^{2}}
\mathbf{1}^{\top}\mu_{22}\mathbf{1}
- \frac{1}{8\sigma^{2}}\sum_{i}\partial_{iiii}C(0)
+o(1)+O(u^{-2}|z|^{2}|y|^{2})
+O(u^{-1}\mathcal{E}^{2}({\mathbf z})).
\end{eqnarray*}
We write $X_u = o_{p}(1)$ if $X_u \rightarrow0$ in probability as
$u\rightarrow\infty$. We insert the above form back to the integral
in (\ref{ProbS}) and apply Lemma \ref{LemRemainder},
\begin{eqnarray*}
p(\Xi_{\varepsilon})
&=&o(1)u^{-\alpha}e^{-u^{2}/2}\\[-2pt]
&&{}+\frac{1}{(2\pi
)^{(d+1)(d+2)/4}}|\Gamma|^{-1/2}\\[-2pt]
&&\hspace*{10pt}{}\times\int_{\mathcal{L}}P\bigl(u\cdot
A>o_{p}(1)\bigr)
\\[-2pt]
&&\hspace*{-24.5pt}\hspace*{34.3pt}{}\times\exp\biggl\{-\biggl[\frac{1}{2}u^{2}+\frac{u}{\sigma}A+\frac
{1}{2}B^{\top}B\\[-2pt]
&&\hspace*{-24.5pt}\hspace*{81.7pt}{}+
\frac{1}{2}\frac{(-{A}/{\sigma}+\mu_{20}\mu_{22}^{-
{1/2}}B+%
{1}/({2\sigma})\mu_{20}\mathbf{1}+o(1))^{2}}
{1-\mu_{20}\mu
_{22}^{-1}\mu_{02}} \\[-2pt]
&&\hspace*{-24.5pt}\hspace*{81.7pt}{} +\frac{u}{\sigma}H_{0}-\frac{u}{\sigma}H\bigl((uI-\mathbf
{z})^{-1/2}y,(uI-%
\mathbf{z})^{1/2},\varepsilon\bigr) \\[-2pt]
&&\hspace*{-24.5pt}\hspace*{81.7pt}{}-\frac{1}{8\sigma^{2}}\mathbf{1}^{\top}\mu_{22}\mathbf
{1}-%
\frac{1}{8\sigma^{2}}\sum_{i}\partial
_{iiii}C(0)+o(1)\\[-2pt]
&&\hspace*{-24.5pt}\hspace*{81.7pt}\hspace*{51.4pt}{}+O(u^{-2}|z|^{2}|y|^{2})+O(u^{-1}\mathcal{E}^{2}(%
{\mathbf z}))\biggr]\biggr\}\,dw\,dy\,dz.
\end{eqnarray*}
Note that Jacobian determinant is
\[
\biggl\vert\det\biggl( \frac{\partial(w,z,y)}{\partial
(A,B,y)}\biggr)
\biggr\vert=\sigma^{-1}\det(\mu_{22})^{1/2}.
\]
Note that when $|(uI-{\mathbf z})^{-1/2}y|\leq\kappa u^{\delta
}-u^{\delta/2}$ (the first situation in Lemma~\ref{LemJ}),
\[
uH_{0}-uH\bigl((uI-\mathbf{z})^{-1/2}y,(uI-\mathbf{z})^{1/2},\varepsilon\bigr)=o(1).
\]
Then, with another change of variable, $A^{\prime}=uA$, the
integration on $%
\mathcal{L}_{1}$ is
%
%
\begin{eqnarray}\label{L1}
&&\int_{\mathcal{L}_{1}}h(w,y,z)P\biggl( \int_{\Xi_{\varepsilon
}}e^{\sigma E(t)+\sigma g(t)}\,dt>b\biggr) \,dw\,dy\,dz \nonumber\\[-2pt]
&&\qquad=\frac{|\Gamma|^{-1/2}}{(2\pi)^{(d+1)(d+2)/4}}\nonumber\\[-2pt]
&&\qquad\quad{}\times\int_{\mathcal{L}%
_{1}}P\bigl( u\cdot A>o_{p}(1)\bigr) \nonumber\\[-2pt]
&&\qquad\quad{}\times\exp\biggl\{-\biggl[\frac{1}{2}u^{2}+\frac{u}{\sigma}A+\frac
{1}{2}B^{\top}B\nonumber\\[-2pt]
&&\qquad\quad\hspace*{-28.6pt}\hspace*{74.88pt}{}+
\frac{1}{2}\frac{(-{A}/{\sigma}+\mu_{20}\mu_{22}^{-
{1/2}}B+%
{1}/({2\sigma})\mu_{20}\mathbf{1}+o(1))^{2}}
{{1-\mu_{20}\mu
_{22}^{-1}\mu_{02}}} \\[-2pt]
&&\qquad\quad\hspace*{-28.6pt}\hspace*{73.8pt}{}+\frac{u}{\sigma}H_{0} -\frac{u}{\sigma}H\bigl((uI-\mathbf
{z})^{-1/2}y,(uI-%
\mathbf{z})^{1/2},\varepsilon\bigr)\nonumber\\[-2pt]
&&\qquad\quad\hspace*{-28.6pt}\hspace*{73.8pt}{}-\frac{1}{8\sigma
^{2}}\mathbf{1}%
^{\top}\mu_{22}\mathbf{1}-\frac{1}{8\sigma^{2}}\sum_{i}\partial
_{iiii}C(0)+o(1) \nonumber\\[-2pt]
&&\qquad\quad\hspace*{-28.6pt}\hspace*{130.6pt}{} +o(u^{-2}|z|^{2}|y|^{2})+O(u^{-1}\mathcal{E}^{2}({\mathbf z}))\biggr]%
\biggr\}\,dw\,dy\,dz \nonumber\\[-2pt]
&&\qquad=\frac{\sigma^{-1}|\Gamma|^{-1/2}}{(2\pi)^{(d+1)(d+2)/4}}\det
(\mu_{22})^{1/2}\nonumber\\[-2pt]
&&\qquad\quad{}\times e^{-%
({1/2})u^{2}+({1}/({8\sigma^{2}}))\mathbf{1}^{\top}\mu
_{22}\mathbf{1}+%
({1}/({8\sigma^{2}}))\sum_{i}\partial_{iiii}C(0)}\nonumber\\[-2pt]
&&\qquad\quad{}\times\int_{\mathcal{L}
_{1}}P\bigl(u\cdot A>o_{p}(1)\bigr) \nonumber\\[-2pt]
&&\qquad\quad{}\times\exp\biggl\{-\biggl[\frac{u}{\sigma}\cdot A+\frac{B^{\top
}B}{2}\nonumber\\[-2pt]
&&\qquad\quad\hspace*{46.7pt}{} +\frac{(-%
{A}/{\sigma}+\mu_{20}\mu_{22}^{-{1/2}}B+
{1}/({2\sigma})\mu
_{20}\mathbf{1}+o(1))^2}
{2(1-\mu_{20}\mu_{22}^{-1}\mu_{02})}
\nonumber\\[-2pt]
&&\qquad\quad\hspace*{46.7pt}\hspace*{30.5pt}{}
+o(1)+O\bigl(u^{-2}|z|^{2}|y|^{2}+u^{-1}\mathcal{E}^{2}({\mathbf z})\bigr)%
\biggr]\biggr\}\,dA\,dB\,dy \nonumber\\[-2pt]
&&\qquad=\frac{\sigma^{-1}|\Gamma|^{-1/2}}{(2\pi)^{(d+1)(d+2)/4}}\det
(\mu
_{22})^{1/2}u^{-1}\nonumber\\[-2pt]
&&\qquad\quad{}\times e^{-({1/2})u^{2}+({1}/({8\sigma^{2}}))\mathbf
{1}%
^{\top}\mu_{22}\mathbf{1}+({1}/({8\sigma^{2}}))\sum_{i}\partial
_{iiii}C(0)}\nonumber\\[-2pt]
&&\qquad\quad{}\times\int_{\mathcal{L}_{1}}P\bigl(A^{\prime}>o_{p}(1)\bigr) \nonumber\\[-2pt]
&&\qquad\quad{}\times\exp\biggl\{-\biggl[\frac{A^{\prime}}{\sigma}+\frac{B^{\top
}B}{2}\nonumber\\[-2pt]
&&\qquad\quad\hspace*{46.9pt}{}+\frac{(-%
{A^{\prime}}/({\sigma u})+\mu_{20}\mu_{22}^{-{1/2}}B+
{1}/({%
2\sigma})\mu_{20}\mathbf{1}+o(1))^{2}}{2(1-\mu_{20}\mu_{22}^{-1}\mu
_{02})%
}\nonumber\\[-2pt]
&&\qquad\quad\hspace*{46.9pt}\hspace*{34.6pt}{} +o(1)+O\bigl(u^{-2}|z|^{2}|y|^{2}+u^{-1}\mathcal{E}^{2}({\mathbf z})\bigr)%
\biggr]\biggr\}\,dA^{\prime}\,dB\,dy. \nonumber
\end{eqnarray}
The second equality is a change of variable from $(w,y,z)$ to
$(A,B,y)$. The
third equality is a change of variable from $(A,B,y)$ to $(A^{\prime
},B,y)$%
. Note that $P(A^{\prime}>o_{p}(1))\rightarrow I(A^{\prime}>0)$ as $%
u\rightarrow\infty$. In addition, on the set $\mathcal{L}$,
\[
O\bigl(u^{-2}|z|^{2}|y|^{2}+u^{-1}\mathcal{E}^{2}({\mathbf
z})\bigr)=O(u^{-1+2\delta+2\varepsilon_{0}}|z|^{2}).
\]
By choosing $\delta$ and $\varepsilon_{0}$ small enough, when
$|B|<u^{1/4}$%
, $u^{-1+2\delta+2\varepsilon_{0}}|z|^{2}=o(1)$; $|B|>u^{1/4}$, $%
|B|=\Theta(|z|)$, therefore,
\[
\frac{B^{\top}B}{2}+O(u^{-1+2\delta+2\varepsilon
_{0}}|z|^{2})=\bigl(1+o(1)\bigr)%
\frac{B^{\top}B}{2}.
\]
The integrant in (\ref{L1}) has the following bound, for $A'>0$
\begin{eqnarray*}
\hspace*{-5pt}&&P\bigl(A^{\prime}>o_{p}(1)\bigr)\\[-2pt]
\hspace*{-5pt}&&\hspace*{-1pt}\quad{}\times
\exp\biggl\{-\biggl[\frac{A^{\prime}}{\sigma}+\frac{%
B^{\top}B}{2}
+ \frac{(-{A^{\prime}}/({\sigma u})+\mu_{20}\mu
_{22}^{-%
{1/2}}B+({1}/({2\sigma}))\mu_{20}\mathbf
{1}+o(1))^{2}}{2(1-\mu
_{20}\mu_{22}^{-1}\mu_{02})} \\[-2pt]
\hspace*{-5pt}&&\hspace*{-1pt}\qquad\quad\hspace*{104.3pt}{} +o(1)+O\bigl(u^{-2}|z|^{2}|y|^{2}+u^{-1}\mathcal{E}^{2}({\mathbf z})\bigr)%
\biggr]\biggr\} \\[-2pt]
\hspace*{-5pt}&&\hspace*{-1pt}\qquad\leq 2\exp\biggl\{\!-\frac{1}{\delta^{\prime}}\biggl[\frac{A^{\prime
}}{%
\sigma}+\frac{B^{\top}B}{2}+\frac{(-{A^{\prime}}/({\sigma
u})+\mu
_{20}\mu_{22}^{-{1/2}}B+({1}/({2\sigma}))\mu_{20}\mathbf
{1})^{2}}{%
2(1-\mu_{20}\mu_{22}^{-1}\mu_{02})}\biggr]\biggr\}
\end{eqnarray*}
for $\delta^{\prime}$ small enough. Note that the $o_p(1)$ is in fact
$-u\log E\exp(g((uI-\zz)^{-1/2}S))$. Thanks to the result of Lemma
\ref{LemRemainder}, the integral of the left-hand side of the above display in
the region $A'<0$ is $o(1)$. By dominated convergence
theorem, (\ref{L1})~equals
%
%
\begin{eqnarray}\label{ApxL1}\quad
&&\bigl(1+o(1)\bigr)\frac{\sigma^{-1}|\Gamma|^{-1/2}}{(2\pi
)^{(d+1)(d+2)/4}}\nonumber\\[-2pt]
&&\quad{}\times\det(\mu
_{22})^{1/2}u^{-1}e^{-({1/2})u^{2}+({1}/({8\sigma^{2}}))\mathbf
{1}%
^{\top}\mu_{22}\mathbf{1}+({1}/({8\sigma^{2}}))\sum_{i}\partial
_{iiii}C(0)} \nonumber\\[-2pt]
&&\quad{}\times\int_{A^{\prime}>0,|y|_{\infty}<\kappa u^{\delta+1/2}-u^{\delta
/2+1/2}}\exp\biggl\{-\biggl[A^{\prime}/\sigma+\frac{B^{\top
}B}{2}\\
&&\qquad\quad\hspace*{97pt}{}+\frac{(\mu
_{20}\mu_{22}^{-{1/2}}B+({1}/({2\sigma}))\mu_{20}\mathbf
{1})^{2}}{%
2(1-\mu_{20}\mu_{22}^{-1}\mu_{02})}\biggr]\biggr\}\,dA^{\prime}\,dB\,dy\hspace*{-1pt}\nonumber\\[-2pt]
&&\qquad=\bigl(1+o(1)\bigr)H\mes(u\Xi_{\varepsilon})u^{-1}e^{-({1/2})u^{2}},
\nonumber
\end{eqnarray}
where $H$ is defined in (\ref{H}). The above display is obtained by
the fact that $\mes(u\Xi_{\varepsilon})= u^{d}\mes(\Xi_{\varepsilon
})=\kappa^du^{d/2+d\delta}$.

\subsection*{\texorpdfstring{Situations 2 and 3 of Lemma \protect\ref{LemJ}}{Situations 2 and 3 of Lemma 5}}

For the second situation in Lemma \ref{LemJ}, let $\mathcal
{L}_{2}=\mathcal{L%
}\cap\{\kappa u^{\delta}-u^{\delta/2}<|(uI-{\mathbf z}%
)^{-1/2}y|_\infty\leq(1+\varepsilon_{1})\kappa u^{\delta}\}$ and
there exists $%
c_{1}>0$ such that
%
%
\begin{eqnarray} \label{ApxL2}
&&\int_{\mathcal{L}_{2}}P\bigl( u\cdot A>o_{p}(1)\bigr)\nonumber\\[-2pt]
&&\quad{}\times \exp\biggl\{
-\biggl[\frac{1}{2}%
u^{2}+\frac{u}{\sigma}A+\frac{1}{2}B^{\top}B\nonumber\\[-2pt]
&&\quad\hspace*{47.6pt}{}+\frac{1}{2}\frac
{(-{A/%
\sigma}+\mu_{20}\mu_{22}^{-{1/2}}B+({1}/({2\sigma}))\mu
_{20}%
\mathbf{1}+o(1))^{2}}{1-\mu_{20}\mu_{22}^{-1}\mu_{02}} \nonumber\\[-2pt]
&&\quad\hspace*{47.6pt}{} +\frac{u}{\sigma}H_{0}-\frac{u}{\sigma}H\bigl((uI-\mathbf
{z})^{-1/2}y,(uI-%
\mathbf{z})^{1/2},\varepsilon\bigr)\\[-2pt]
&&\quad\hspace*{47.6pt}{}-\frac{1}{8\sigma
^{2}}\mathbf{1}%
^{\top}\mu_{22}\mathbf{1}-\frac{1}{8\sigma^{2}}\sum_{i}\partial
_{iiii}C(0)+o(1) \nonumber\\[-2pt]
&&\quad\hspace*{134.7pt}{} +o(u^{-2}|z|^{2}|y|^{2})+O(u^{-1}\mathcal{E}^{2}({\mathbf z}%
))\biggr]\biggr\}\,dw\,dy\,dz \nonumber\\[-2pt]
&&\qquad\leq \bigl(c_{1}\varepsilon_{1}^{d}+o(1)\bigr)H \mes(u\Xi_{\varepsilon
})u^{-1}e^{-%
({1/2})u^{2}}.\nonumber\vadjust{\goodbreak}
\end{eqnarray}
For the third situation, $\mathcal{L}_{3}=\mathcal{L}\cap\{
(1+\varepsilon
_{1})\kappa u^{\delta}<|(uI-{\mathbf z})^{-1/2}y|\leq
u^{\delta
+\varepsilon_{0}}\}$,
%
%
\begin{eqnarray}\label{ApxL3}
&&\int_{\mathcal{L}_{3}}P\bigl( u\cdot A>o_{p}(1)\bigr)\nonumber\\
&&\quad{}\times \exp\biggl\{
-\biggl[\frac{1}{2}%
u^{2}+\frac{u}{\sigma}A+\frac{1}{2}B^{\top}B\nonumber\\
&&\quad\hspace*{47.3pt}{}+\frac{1}{2}\frac
{(-{A/%
\sigma}+\mu_{20}\mu_{22}^{-{1/2}}B+({1}/({2\sigma}))\mu
_{20}%
\mathbf{1}+o(1))^{2}}{1-\mu_{20}\mu_{22}^{-1}\mu_{02}} \nonumber\\
&&\quad\hspace*{47.3pt}{} +\frac{u}{\sigma}H_{0}-\frac{u}{\sigma}H\bigl((uI-\mathbf
{z})^{-1/2}y,(uI-%
\mathbf{z})^{1/2},\varepsilon\bigr)\nonumber\\[-8pt]\\[-8pt]
&&\quad\hspace*{47.3pt}{}-\frac{1}{8\sigma
^{2}}\mathbf{1}%
^{\top}\mu_{22}\mathbf{1}-\frac{1}{8\sigma^{2}}\sum_{i}\partial
_{iiii}C(0)+o(1) \nonumber\\
&&\quad\hspace*{135pt}{} +o(u^{-2}|z|^{2}|y|^{2})+O(u^{-1}\mathcal{E}^{2}({\mathbf z}%
))\biggr]\biggr\}\,dw\,dy\,dz \nonumber\\
&&\qquad\leq O(1)\biggl(\frac{1}{2}\biggr)^{u/\sigma}u^{-1}u^{({1/2}+\delta
+\varepsilon
_{0})d}e^{-({1/2})u^{2}}\nonumber\\
&&\qquad=o(1)\mes(u\Xi_{\varepsilon
})u^{-1}e^{-(1/2)u^{2}}.\nonumber
\end{eqnarray}
We put (\ref{ApxL1}), (\ref{ApxL2}) and (\ref{ApxL3}) together and conclude
the proof.\
\end{pf*}

\section{\texorpdfstring{Proof for Theorem \protect\ref{PropInd}}{Proof for Theorem 3}}\label{sec_proof2}

Similar to Section \ref{sec_proof1}, we arrange
all the lemmas and their proofs in the \hyperref[app]{Appendix}.
\begin{pf*}{Proof of Theorem \ref{PropInd}}
Since the proofs for $\mathcal C^+$ and $\mathcal C^-$ are
complete analogue, we only provide the proof for $\mathcal
C^+$. We prove for the asymptotics by providing bounds from
both sides. We first discuss the easy case: the lower bound.
Note that
\begin{eqnarray*}
&&
P\biggl( I_\sigma\biggl(\bigcup_{\kk\in\mathcal{C}^{+}}\Xi
_{\varepsilon,\kk}\biggr)>b\biggr)\\
&&\qquad\geq P\Bigl(\max_{\kk\in\mathcal
C^{+}}I_\sigma(\Xi_{\varepsilon,\kk})>b\Bigr) \\
&&\qquad\geq
\sum_{\kk\in\mathcal C^{+}}P\bigl(I_\sigma(\Xi_{\varepsilon
,\kk})>b\bigr)-\sum_{\kk\neq\kk'}P\bigl(I_\sigma(\Xi_{\varepsilon
,\kk})>b,I_\sigma(\Xi_{\varepsilon,\kk'})>b\bigr).
\end{eqnarray*}
Thanks to Lemma \ref{LemPair},
\[
P\biggl( I_\sigma\biggl(\bigcup_{\kk\in
\mathcal{C}^{+}}\Xi_{\varepsilon,\kk}\biggr)>b\biggr)\geq
\bigl(1+o(1)\bigr)\sum_{\kk\in\mathcal C^{+}}P\bigl(I_\sigma(\Xi
_{\varepsilon,\kk})>b\bigr).
\]
The rest of the proof
is to establish the asymptotic upper bound. To simplify our
writing, we let
%
%
\begin{equation}\label{DefABD}
\cA= I_\sigma(\Xi_\varepsilon),\qquad
\cB=
I_\sigma\Bigl(\sup_{\kk\in\mathcal N}\Xi_{\varepsilon,\kk}\Bigr),\qquad
\cD= I_\sigma\biggl(\bigcup_{\kk'\in\cC^+ \setminus[\{ \mathbf
0\}\cup\mathcal N]
}\Xi_{\varepsilon,\kk'}\biggr),\hspace*{-35pt}
\end{equation}
where
$\mathcal N$ is the set of neighbors of $\Xi_\varepsilon$,
that is, $\kk\in\mathcal N$ if and only if
\[
{\inf_{s\in
\Xi_\varepsilon, t\in\Xi_{\varepsilon, \kk} }}|s-t|=0.
\]
An illustration of $\cA$, $\cB$ and $\cD$ is given in Figure
\ref{FigSQ}.

%
\begin{figure}

\includegraphics{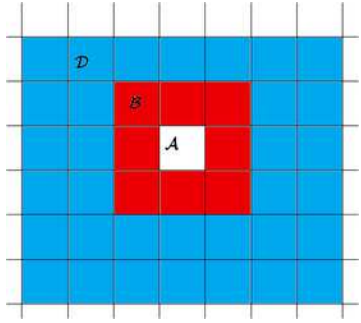}

\caption{Illustration of $\mathcal A$, $\mathcal B$ and
$\mathcal D$.} \label{FigSQ}
\vspace*{6pt}
\end{figure}

Further, let
%
%
\begin{equation}\label{B0}
b_0 = u^{-1-d/2}b, \qquad
b-b_0 = (1-u^{-1-d/2})b = e^{-(1+o(1))u^{-1-d/2}}b,\hspace*{-35pt}
\end{equation}
and
$u_0$ solves
\[
u_0^{-d/2}e^{\sigma u_0+H_{0}}=b_0,
\]
and there exists $c_0>0$ such that $u_0>u-c_0\log u$. The first
step in developing the upper bound is to use the following
inequality
%
%
\begin{eqnarray}\label{upper}
&&
P(\cA+ \cB+\cD>b)\nonumber\\[2pt]
&&\qquad\leq P(\cA> b-b_0) + P(\cA\leq
b_0, \cA+ \cB+\cD>b)\nonumber\\[2pt]
&&\qquad\quad{} + P(b_0 < \cA\leq b-b_0, \cA+ \cB+\cD
>b)\\[2pt]
&&\qquad\leq P(\cA> b-b_0)+ P(\cB+\cD>b-b_0) \nonumber\\[2pt]
&&\qquad\quad{}+ P(\cA>b_0,\cB+
\cD>b_0, \cA+ \cB+ \cD
> b).\nonumber
\end{eqnarray}
From Theorem \ref{PropSmall},
\[
P(\cA> b-b_0) =
\bigl(1+o(1)\bigr)P(\cA> b).
\]
The next step is to show that the last
term in (\ref{upper}) is ignorable. Note that
\begin{eqnarray*}
&&P( \cA>b_{0},\cB+\cD>b_{0},\cA+\cB+\cD>b) \\[-2pt]
&&\qquad=P( \cA+\cB>b-b_{0},\cA>b_{0},\cB+\cD>b_{0},\cA+\cB+\cD
>b) \\[-2pt]
&&\qquad\quad{} +P( \cD>b-b_{0},\cA>b_{0},\cB+\cD>b_{0},\cA+\cB+\cD
>b) \\[-2pt]
&&\qquad\quad{}+P( \cA+\cB>b_{0},\cD>b_{0},\cA>b_{0},\cB+\cD>b_{0},\cA
+\cB+\cD>b)  \\[-2pt]
&&\qquad\leq P(
\cA+\cB>b-b_{0},\cA>b_{0},\cB+\cD>b_{0},\cA+\cB+\cD>b)\\[-2pt]
&&\qquad\quad{} + 2P( \cD>b_{0},\cA>b_{0})\\[-2pt]
&&\qquad= o\bigl(P( \cA>b)\bigr).
\end{eqnarray*}
The last step is due to Lemmas \ref{LemAD} and \ref{Lemlast}.
By noting that $\#(\cC^+) u^{-1-d/2} = o(u^{-1})$, the
conclusion of the theorem is immediate by induction, where
$\#(\cdot)$ is the cardinality of a set.\vspace*{-3pt}
\end{pf*}

\begin{appendix}\label{app}

\section*{\texorpdfstring{Appendix: Lemmas in Sections \lowercase{\protect\ref{sec_proof1}} and
\lowercase{\protect\ref{sec_proof2}}}{Appendix: Lemmas in Sections 4 and 5}}

Lemma \ref{LemEx1} isolated the dominating event so that we will be
in good shape to use Taylor's expansion.
\begin{lemma}\label{LemEx1}
There exist $\varepsilon_0, \delta>0$ small enough and
$\kappa$ large. Let $\varepsilon= \kappa u^{-1/2
+\delta}$ such that for any $\alpha
>0$,
\begin{eqnarray*}
&&P\biggl( \vert f( 0) -u\vert>u^{2\delta
+\varepsilon_{0}}\mbox{ or }\vert\partial f( 0)
\vert
>u^{{1}/{2}+\delta+\varepsilon_{0}}\mbox{ or }\\[-2pt]
&&\hspace*{21.2pt}\vert
\partial
^{2}f(0)-u\mu_{20}\vert>u^{{1}/{2}+\varepsilon
_{0}},\int{\Xi_\varepsilon}e^{f( t) }\,dt>b\biggr) \\[-2pt]
&&\qquad=o( 1) u^{-\alpha}e^{-({1}/{2})u^{2}}.
\end{eqnarray*}
\end{lemma}
\begin{pf} 
Note that there exists $c_{1}$ such that $\sigma u\leq\log
b+c_{1}\log\log b$. Let $\sigma\tilde{u}=\log(b)$. Since we
only consider the case that $u$ is large, we always have
\mbox{$\mes(\Xi_{\varepsilon})<1$}:
\begin{eqnarray*}
&&P\biggl( f(0)<u-u^{2\delta+\varepsilon_{0}},\int_{\Xi
_{\varepsilon
}}e^{f(t)}\,dt>b\biggr) \\[-2pt]
&&\qquad\leq P\bigl(f(0)<u-u^{2\delta+\varepsilon_{0}},\sup f(t)>\tilde{u}\bigr) \\[-2pt]
&&\qquad\leq CP\bigl(f(0)<u-u^{2\delta+\varepsilon_{0}}|\sup f(t)>\tilde
{u}\bigr)\tilde{u}%
^{d-1}e^{-\tilde{u}^{2}/2}.
\end{eqnarray*}
The last inequality is an application of Proposition
\ref{ThmPiterbarg}. Because for any $u^{\prime}>\tilde{u}$,
for some $\varepsilon_{1}>0$,
\begin{eqnarray*}
&&\inf_{t\in\Xi_{\varepsilon
}}E\Bigl(f(t)\big|\sup_{\Xi_\varepsilon} f(t)=u^{\prime}\Bigr)\\[-2pt]
&&\qquad\geq
u^{\prime}\inf_{t\in\Xi_{\varepsilon}}C(t)\geq u-\kappa
^{2}\varepsilon_{1}u^{2\delta}\bigl(1+o(1)\bigr)\vadjust{\goodbreak}
\end{eqnarray*}
and
\[
\sup_{t\in(-\varepsilon,\varepsilon
)}\operatorname{Var}\Bigl(f(t)\big|\sup_{\Xi_\varepsilon} f(t)=u^{\prime
}\Bigr)=O(\varepsilon^{2})=O(u^{-1+2\delta}),
\]
one can choose $\kappa$ large enough such that
\[
P\Bigl(f(0)<u-u^{2\delta+\varepsilon_{0}}\big|\sup_{\Xi_\varepsilon}
f(t)>\tilde{u}\Bigr)=O(1)\exp(-{%
u^{1+\varepsilon_{0}}}/{c_{2}}).
\]
Therefore,
\[
P\biggl(f(0)<u-u^{2\delta+\varepsilon_{0}},\int_{\Xi
_{\varepsilon}}e^{f(t)}\,dt>b\biggr)=o(1)u^{-\alpha}e^{-u^{2}/2}
\]
for all $\alpha>0$. Also, $P(f(0)>u+u^{2\delta+\varepsilon
_{0}})=o(1)u^{-\alpha}e^{-({1}/{2})u^{2}}$. Hence,%
\[
P\biggl(|f(0)-u|>u^{2\delta+\varepsilon_{0}},\int_{\Xi
_{\varepsilon}}e^{f(t)}\,dt>b\biggr)=o(1)u^{-\alpha
}e^{-({1}/{2})u^{2}}.
\]
Similarly, we have
\begin{eqnarray*}
&&P\biggl( |f(0)-u|<u^{2\delta+\varepsilon_{0}},|\partial
f(0)|>u^{{1}/{2}+\delta+\varepsilon_{0}},\int_{\Xi
_{\varepsilon}}e^{f(t)}\,dt>b\biggr)\\
&&\qquad=o(1)u^{-\alpha}e^{-({1}/{2})u^{2}}, \\
&&P\biggl( |f(0)-u|<u^{2\delta+\varepsilon_{0}},|\partial
^{2}f(0)-u\mu_{20}|>u^{{1}/{2}+\varepsilon
_{0}},\int_{\Xi_{\varepsilon}}e^{f(t)}\,dt>b\biggr)\\
&&\qquad= o(1)u^{-\alpha}e^{-({1}/{2})u^{2}}.
\end{eqnarray*}
The above two displays are immediate by noting that $
(f(0),\partial f(0),\partial^2 f(0))$ is a multivariate
Gaussian random vector. In addition, $ (f(0),\partial^2 f(0))$
is independent of $\partial f(0)$ and the covariance between
$f(0)$ and $\partial^2 f(0)$ is $\mu_{02}$.
\end{pf}
\begin{lemma}\label{LemDensity}
Let $h(w,y,z)$ be the density of $(f(0), \partial f(0),
\partial^2 f(0))$. Then,
\begin{eqnarray*}
\hspace*{-4pt}&&h(w,y,z)\\
\hspace*{-4pt}&&\qquad=\frac{1}{(2\pi)^{(d+1)(d+2)/4}}|\Gamma
|^{-1/2}\\
\hspace*{-4pt}&&\qquad\quad{}\times\exp\biggl\{-\biggl[\frac{1
}{2}u^{2}+\frac u\sigma A (w,y,z)+\frac{1}{2}y^\top\bigl(I-(I-\zz
/u)^{-1}\bigr) y \\
\hspace*{-4pt}&&\qquad\quad\hspace*{47.6pt}{} +\frac{1}{2}\frac{(w+\mu_{20}\mu_{22}^{-1}z)^{2}}{1-\mu
_{20}\mu_{22}^{-1}\mu_{02}}+\frac{1}{2}z^{\top}\mu
_{22}^{-1}z+\frac{u}{2\sigma}\log
\det(I-u^{-1}\zz) \\
\hspace*{-4pt}&&\qquad\quad\hspace*{121pt}{} -\frac u\sigma\log\int_{|(uI-\zz)^{-1/2}t|<\varepsilon
}e^{\sigma J(t)}\,dt+uH_{0}/\sigma\biggr]\biggr\}.
\end{eqnarray*}
In addition,
%
%
\begin{equation} \label{gammainv}
\Gamma^{-1}=\pmatrix{
\dfrac{1}{1-\mu_{20}\mu_{22}^{-1}\mu_{02}} & -\dfrac{\mu
_{20}\mu_{22}^{-1}
}{1-\mu_{20}\mu_{22}^{-1}\mu_{02}} \vspace*{3pt}\cr
-\dfrac{\mu_{22}^{-1}\mu_{02}}{1-\mu_{20}\mu
_{22}^{-1}\mu_{02}} & \mu_{22}^{-1}+\dfrac{\mu_{22}^{-1}\mu
_{02}\mu_{20}\mu_{22}^{-1}}{1-\mu
_{20}\mu_{22}^{-1}\mu_{02}}}.
\end{equation}
\end{lemma}
\begin{pf}
The form of $\Gamma^{-1}$ in (\ref{gammainv}) is a result from
linear algebra. The form of $\Gamma^{-1}$ is direct
application of the block matrix inverse from linear algebra.
Note that
\begin{eqnarray*}
h(w,y,z)&=&\frac{1}{(2\pi)^{(d+1)(d+2)/4}}|\Gamma|^{-1/2}\\
&&{}\times\exp
\left\{ -\frac{1}{2}(u-w,z^{\top}+u\mu_{20},y^{\top})
\pmatrix{\Gamma^{-1} & 0 \cr
0 & I}
\pmatrix{u-w \cr
z+u\mu_{02} \cr
y}\right\}.
\end{eqnarray*}
By plugging in the form of $\Gamma^{-1}$, we have
\begin{eqnarray*}
&&(u-w,z^{\top}+u\mu_{20},y^{\top})
\pmatrix{\Gamma^{-1} & 0 \cr
0 & I}\pmatrix{u-w \cr
z+u\mu_{02} \cr y}
\\
&&\qquad = y^{\top}y+\frac{(u-w)^{2}}{1-\mu_{20}\mu_{22}^{-1}\mu
_{02}}\\
&&\qquad\quad{}+(z+u\mu_{02})^{\top}\biggl[\mu_{22}^{-1}+\frac{\mu
_{22}^{-1}\mu_{02}\mu_{20}\mu
_{22}^{-1}}{1-\mu_{20}\mu_{22}^{-1}\mu_{02}}\biggr](z+u\mu_{02}) \\
&&\qquad\quad{} -2(u-w)\frac{\mu_{20}\mu_{22}^{-1}}{1-\mu_{20}\mu
_{22}^{-1}\mu_{02}}
(z+u\mu_{02}) \\
&&\qquad= u^{2}+y^{\top}y-2wu+\frac{w^{2}}{1-\mu_{20}\mu_{22}^{-1}\mu
_{02}} \\
&&\qquad\quad{} +z^{\top}\biggl[\mu_{22}^{-1}+\frac{\mu_{22}^{-1}\mu_{02}\mu
_{20}\mu_{22}^{-1}}{1-\mu_{20}\mu_{22}^{-1}\mu
_{02}}\biggr]z+2w\frac{\mu_{20}\mu
_{22}^{-1}}{1-\mu_{20}\mu_{22}^{-1}\mu_{02}}z \\
&&\qquad= u^{2}+\frac{2u}{\sigma}A(w,y,z)+y^{\top
}\bigl(I-(uI-\zz)^{-1}\bigr)y+\frac{(w+\mu_{20}\mu
_{22}^{-1}z)^{2}}{1-\mu_{20}\mu_{22}^{-1}\mu_{02}} \\
&&\qquad\quad{} +z^{\top}\mu_{22}^{-1}z+\frac u\sigma\log\det(I-u^{-1}\zz)\\
&&\qquad\quad{}-\frac{2u}{\sigma}\log
\int_{|(uI-\zz)^{-1/2}t|<\varepsilon}e^{J(t)}\,dt+\frac{2u}{\sigma}H_{0}.
\end{eqnarray*}
Therefore, we conclude the proof.\vadjust{\goodbreak}
\end{pf}
\begin{lemma}\label{LemDet}
\[
\log\bigl(\det(I-u^{-1}\zz)\bigr) = - u^{-1}\operatorname{Tr}(\zz) + \tfrac1 2
u^{-2}\mathcal E^2(\zz) + o(u^{-2}),
\]
where $\mathrm{Tr}$ is the trace
of a matrix, $\mathcal E^2 (\zz)=\sum_{i=1}^d \lambda_i^2$, and
$\lambda_i$'s are the eigenvalues of~$\zz$.
\end{lemma}
\begin{pf}
The result is immediate by noting that
\[
\det(I-u^{-1}\zz)=\prod_{i=1}^{d}(1-\lambda_{i}/u)
\]
and $\operatorname{Tr}(\zz)= \sum_{i=1}^d \lambda_i$.
\end{pf}
\begin{lemma}\label{lemG}
Let $y=(y_1,\ldots,y_d)^\top$ and $X\sim N(y/\sqrt u,I/\sqrt
\sigma)$. Then, on the set $\mathcal L$ defined in (\ref{L}),
\begin{eqnarray*}
&&E\bigl(g_{3}\bigl(X/\sqrt{u}\bigr)+g_{4}\bigl(X/\sqrt{u}\bigr)\bigr)\\
&&\qquad=-\frac
{u^{-1}}{8}(u^{-1}Y+\mathbf{1}/\sigma)
^{\top}\mu
_{22}(u^{-1}Y+\mathbf{1}/\sigma)\\
&&\qquad\quad{}
+\frac{u^{-1}}{8\sigma^2}\mathbf{1}^{\top
}\mu_{22}\mathbf{1}+\frac{u^{-1}}{8\sigma^2}\sum_{i}\partial
_{iiii}C(0) +o(u^{-1}),
\end{eqnarray*}
where
\begin{eqnarray*}
Y &=&(y_{i}^{2},i=1,\ldots,d,2y_{i}y_{j},1\leq i<j\leq d)^{\top}, \\
\mathbf{1}&=&(\underbrace{1,\ldots,1}_{d},
\underbrace{0,\ldots,0}_{d(d-1)/2})^{\top}.
\end{eqnarray*}
\end{lemma}
\begin{pf} 
Using the derivatives in (\ref{Derivative}), we
have that
\begin{eqnarray*}
\partial_{ijk}E(0) &=& -\sum_{l=1}^d
\partial^4_{ijkl} C(0)y_l,\\
\partial^4_{ijk}E(0) &=& \bigl(u +
O(|z|+|w|)\bigr)\partial_{ijkl}C(0).
\end{eqnarray*}
We plug this into the
definition of $g_3$ and $g_4$ in (\ref{G34}) and obtain, on the
set $\mathcal L$,
\begin{eqnarray*}
&&E\bigl(g_{3}\bigl(X/\sqrt{u}\bigr)+g_{4}\bigl(X/\sqrt{u}\bigr)\bigr) \\[-2pt]
&&\qquad=-\frac{1}{6}u^{-3/2}\sum_{ijkl}\partial^4
_{ijkl}C(0)E(X_{i}X_{j}X_{k}y_{l})\\[-2pt]
&&\qquad\quad{}+\frac{u^{-1}}{24}\sum
_{ijkl}\partial^4
_{ijkl}C(0)E(X_{i}X_{j}X_{k}X_{l}) +o(u^{-1})\\[-2pt]
&&\qquad=-\frac{1}{8}u^{-3/2}\sum_{ijkl}\partial^4
_{ijkl}C(0)E(X_{i}X_{j}X_{k}y_{l})\\[-2pt]
&&\qquad\quad{}+\frac{u^{-1}}{24}\sum
_{ijkl}\partial^4
_{ijkl}C(0)E\bigl(X_{i}X_{j}X_{k}\bigl(X_{l}-y_{l}/\sqrt{u}\bigr)\bigr) +o(u^{-1})\\[-2pt]
&&\qquad=-\frac{1}{8}u^{-3}\sum_{ijkl}\partial^4
_{ijkl}C(0)y_{i}y_{j}y_{k}y_{l}-
\frac{3}{8\sigma}u^{-2}\sum_{il}y_{i}y_{l}\sum_{j}\partial^4
_{iljj}C(0) \\[-2pt]
&&\qquad\quad{} +\frac{u^{-2}}{8\sigma}\sum_{ij}y_{i}y_{j}\sum_{l}\partial^4
_{ijll}C(0)+\frac{
3u^{-1}}{24\sigma^2}\sum_{i}\partial^4 _{iiii}C(0) +o(u^{-1})\\[-2pt]
&&\qquad=-\frac{1}{8}u^{-3}\sum_{ijkl}\partial^4
_{ijkl}C(0)y_{i}y_{j}y_{k}y_{l}-
\frac{1}{4\sigma}u^{-2}\sum_{ij}y_{i}y_{j}\sum_{l}\partial^4
_{ijll}C(0)\\[-2pt]
&&\qquad\quad{}+\frac{u^{-1} }{8\sigma^2}\sum_{i}\partial^4
_{iiii}C(0) +o(u^{-1}).
\end{eqnarray*}
This last\vspace*{2pt} step is true because
$\sum_{il}y_{i}y_{l}\sum_{j}\partial^4
_{iljj}C(0)=\sum_{ij}y_{i}y_{j}\sum_{l}\partial^4 _{ijll}C(0)$,
which is just a change of index. Then, with the definition of
$Y$ and $\mathbf1$ in the statement of this lemma (note that
$\mathbf1$ is NOT a vector of $1$'s), we have
\begin{eqnarray*}
&&E\bigl(g_{3}\bigl(X/\sqrt{u}\bigr)+g_{4}\bigl(X/\sqrt{u}\bigr)\bigr) \\[-2pt]
&&\qquad=-\frac{u^{-3}}{8}Y^{\top}\mu
_{22}Y-\frac{u^{-2}}{4\sigma}Y^{\top}\mu_{22}
\mathbf{1}+\frac{u^{-1}}{8\sigma^2}\sum_{i}\partial^4 _{iiii}C(0)
+o(1)\\[-2pt]
&&\qquad=-\frac{u^{-1}}{8}(u^{-1}Y+\mathbf{1}/\sigma)^{\top}\mu
_{22}(u^{-1}Y+\mathbf{1}/\sigma)+\frac{u^{-1}}{8\sigma^2}\mathbf
{1}^{\top
}\mu_{22}\mathbf{1}\\[-2pt]
&&\qquad\quad{}+\frac{u^{-1}}{8\sigma^2}
\sum_{i}\partial^4 _{iiii}C(0)+o(1).
\end{eqnarray*}
\upqed\end{pf}
\begin{lemma}\label{LemJ}
Let $J(t)$ be defined in (\ref{J}). Then, on the set $\mathcal
L$ the approximations of $\int_{|(uI-\mathbf
z)^{-1/2}t|_\infty<\varepsilon}e^{\sigma J(t)}\,dt$ under different
situations are as follows:
\begin{longlist}[(2)]
\item[(1)] When $|(uI-\zz)^{-1/2}y|_\infty\leq\kappa u^{\delta
}-u^{\delta/2}$,
\begin{eqnarray*}
&&\int_{|(uI-\mathbf{z})^{-1/2}t|_\infty<\varepsilon}e^{\sigma
J(t)}\,dt \\[-2pt]
&&\qquad= \exp\bigl[
\sigma E\bigl(g_{3}\bigl(X/\sqrt{u}\bigr)+g_{4}\bigl(X/\sqrt{u}\bigr)\bigr)\\[-2pt]
&&\qquad\quad\hspace*{17.6pt}{}+H\bigl((uI-\mathbf{z})^{-1/2}y,
(uI-\mathbf{z})^{1/2},\varepsilon\bigr)+o(u^{-1})\bigr],
\end{eqnarray*}
where
\[
e^{H(y,\Sigma,\varepsilon)}=\int_{|\Sigma^{-1}t|<\varepsilon} e^{
-({\sigma}/{2})(t-y)^{\top}(t-y)}\,dt,\vadjust{\goodbreak}
\]
and $X$ is the random vector defined in Lemma \ref{lemG}.
In addition,
\[
H\bigl(0,(uI-\mathbf{z})^{1/2},\varepsilon\bigr) - H_0 = o(u^{-1}).
\]

\item[(2)] For any $\varepsilon_1 >0$, when $\kappa
u^{\delta}-u^{\delta/2}\leq|(uI-\zz)^{-1/2}y|_\infty\leq
(1+\varepsilon_1)\kappa u^{\delta}$,
\begin{eqnarray*}
&&\int_{|(uI-\mathbf{z})^{-1/2} t|<\varepsilon}e^{\sigma J(t)}\,dt \\[-2pt]
&&\qquad\leq \exp\bigl[
\sigma E\bigl(g_{3}\bigl(X/\sqrt{u}\bigr)+g_{4}\bigl(X/\sqrt{u}\bigr)\bigr)\\[-2pt]
&&\qquad\quad\hspace*{17pt}{} +H\bigl((uI-\mathbf{z})^{-1/2}y,
(uI-\mathbf{z})^{1/2},\varepsilon\bigr)+o(u^{-1})\bigr].
\end{eqnarray*}

\item[(3)] When $(1+\varepsilon_1 )\kappa u^{\delta}<|(uI-\zz
)^{-1/2}y|_\infty\leq u^{\delta+\varepsilon_0}$
\begin{eqnarray*}
&&\int_{|(uI-\mathbf{z})^{-1/2} t|<\varepsilon}e^{\sigma J(t)}\,dt \\[-2pt]
&&\qquad\leq \frac1 2 \exp\bigl[
\sigma E\bigl(g_{3}\bigl(X/\sqrt{u}\bigr)+g_{4}\bigl(X/\sqrt{u}\bigr)\bigr)+H_0+o(u^{-1})\bigr].
\end{eqnarray*}
\end{longlist}
\end{lemma}
\begin{pf}
Note that
\begin{eqnarray*}
&&\int_{|(uI-\mathbf{{\mathbf z}})^{-1/2} t|<\varepsilon
}e^{\sigma J(t)}\,dt \\[-2pt]
&&\qquad=e^{H((uI-{\mathbf z})^{-1/2}y,(uI-{\mathbf z}%
)^{-1/2},\varepsilon)} \\[-2pt]
&&\qquad\quad\hspace*{0pt}{}\times
E\bigl\{ \exp\bigl[ \sigma g_{3}\bigl((uI-{\mathbf z})^{-1/2}X^{\prime}\bigr)\\[-2pt]
&&\qquad\quad\hspace*{24pt}{}+\sigma
g_{4}\bigl((uI-{\mathbf z})^{-1/2}X^{\prime}\bigr)+\sigma R\bigl((uI-{\mathbf z}
)^{-1/2}X^{\prime}\bigr)\bigr] ;\\[-2pt]
&&\qquad\quad\hspace*{149.2pt}|(uI-\mathbf{{\mathbf z}}
)^{-1/2}X^{\prime}|<\varepsilon\bigr\} .
\end{eqnarray*}
Also, $(uI-{\mathbf z})^{-1/2}X'=(1+O(z/u))X^{\prime
}/\sqrt n$ and
\[
X^{\prime}=X-y/\sqrt{u}+(uI-{\mathbf z})^{-1/2}y=X+O(u^{-3/2})|
{\mathbf z}y|.
\]
For the first situation, $|(uI-{\mathbf z})^{-1/2}y|_\infty
\leq\kappa
u^{\delta}-u^{\delta/2}$ and $|z|<u^{1/2+\varepsilon_{0}}$,%
\begin{eqnarray*}
&&E\bigl\{ \exp\bigl[ \sigma g_{3}\bigl((uI-{\mathbf z})^{-1/2}X^{\prime}\bigr)+\sigma
g_{4}\bigl((uI-{\mathbf z})^{-1/2}X^{\prime}\bigr)+\sigma R\bigl((uI-{\mathbf z}
)^{-1/2}X^{\prime}\bigr)\bigr] ;\\[-2pt]
&&\hspace*{252.8pt}|(uI-\mathbf{z}
)^{-1/2}X^{\prime}|<\varepsilon\bigr\} \\[-2pt]
&&\qquad=E\bigl\{ \exp\bigl[ \sigma g_{3}\bigl(X^{\prime
}/\sqrt{u}\bigr)+\sigma g_{4}\bigl(X^{\prime}/\sqrt{u
}\bigr)+\sigma R\bigl(X^{\prime}/\sqrt{u}\bigr)+o(u^{-1})\bigr] ;\\[-2pt]
&&\hspace*{219.2pt}|(uI-\mathbf{z}%
)^{-1/2}X^{\prime}|<\varepsilon\bigr\} \\[-2pt]
&&\qquad=E\bigl\{ \exp\bigl[ \sigma g_{3}\bigl(X^{\prime
}/\sqrt{u}\bigr)+\sigma g_{4}\bigl(X^{\prime}/\sqrt{u
}\bigr)+\sigma R\bigl(X^{\prime}/\sqrt{u}\bigr)+o(u^{-1})\bigr] \bigr\} \\[-2pt]
&&\qquad=E\bigl\{ \exp\bigl[ \sigma g_{3}\bigl(X/\sqrt{u}\bigr)+\sigma
g_{4}\bigl(X/\sqrt{u}\bigr)+\sigma R\bigl(X/\sqrt{u} \bigr)+o(u^{-1})\bigr]
\bigr\} .
\end{eqnarray*}
In addition, because $|R(t)|=O(u^{1/2+\delta}t^{5})$, then%
\begin{eqnarray*}
&&E\bigl\{ \exp\bigl[ \sigma g_{3}\bigl(X/\sqrt{u}\bigr)+\sigma
g_{4}\bigl(X/\sqrt{u}\bigr)+\sigma R\bigl(X/\sqrt{u}\bigr)
\bigr] \bigr\} \\
&&\qquad=E\bigl\{ \exp\bigl[ \sigma g_{3}\bigl(X/\sqrt{u}\bigr)+\sigma
g_{4}\bigl(X/\sqrt{u}\bigr)+o(u^{-1})\bigr] \bigr\}.\vadjust{\goodbreak}
\end{eqnarray*}
Further, $g_{3}(X/\sqrt{u})=o(u^{-1/2+\delta})$ and $g_{4}(X/\sqrt
{u})=o(u^{-1/2+\delta})$, then by repeatedly using Talyor's expansion,
we have
\begin{eqnarray*}
&&E\bigl\{ \exp\bigl[ \sigma g_{3}\bigl(X^{\prime}/\sqrt{u}\bigr)+\sigma
g_{4}\bigl(X^{\prime}/\sqrt{u}%
\bigr)+\sigma R\bigl(X^{\prime}/\sqrt{u}\bigr)\bigr]
;|(uI-\mathbf{z})^{-1/2}X^{\prime
}|<\varepsilon\bigr\} \\
&&\qquad=\exp\bigl\{ E\bigl[ \sigma g_{3}\bigl(X/\sqrt{u}\bigr)+\sigma
g_{4}\bigl(X/\sqrt{u}\bigr)+o(u^{-1})\bigr] \bigr\}.
\end{eqnarray*}
For the second situation, the inequality is immediate by noting
that
\begin{eqnarray*}
&&E\bigl\{ \exp\bigl[ \sigma g_{3}\bigl(X^{\prime
}/\sqrt{u}\bigr)+\sigma g_{4}\bigl(X^{\prime}/\sqrt{u} \bigr)+\sigma
R\bigl(X^{\prime}/\sqrt{u}\bigr)\bigr]
;|(uI-\mathbf{z})^{-1/2}X^{\prime
}|<\varepsilon\bigr\} \\
&&\qquad\leq E\bigl\{ \exp\bigl[ \sigma g_{3}\bigl(X^{\prime
}/\sqrt{u}\bigr)+\sigma g_{4}\bigl(X^{\prime}/ \sqrt{u}\bigr)+\sigma
R\bigl(X^{\prime}/\sqrt{u}\bigr)\bigr] \bigr\}.
\end{eqnarray*}
For the third situation, note that the integral is not focusing
on the dominating part, and the conclusion follows immediately.
\end{pf}

The next lemma is known as the Borel-TIS lemma, which was proved
independently by \cite{Bor75,CIS}.
\begin{lemma}[(Borel-TIS)]\label{LemBorel}
Let $f(t)$, $t\in\mathcal U$, $\mathcal U$ is a parameter set,
be a mean zero Gaussian random field. $f$ is almost surely
bounded on $\mathcal U$. Then,
\[
E\Bigl(\sup_{\mathcal U}f(t) \Bigr)<\infty
\]
and
\[
P\Bigl(\max_{t\in\mathcal{U}}f( t)
-E\Bigl[\max_{t\in\mathcal{U}}f( t) \Bigr]\geq b\Bigr)\leq e^{
-{b^{2}}/({2\sigma_{\mathcal{U}}^{2}})} ,
\]
where
\[
\sigma_{\mathcal{U}}^{2}=\max_{t\in\mathcal{U}}\operatorname{Var}[f( t)].
\]
\end{lemma}
\begin{lemma}\label{LemRemainder}
Let $\log E\exp(g((uI-\zz)^{-1/2}S))$ be defined in
(\ref{Remainder}), then there exists a $\lambda>0$ such that for all $x>0$
\[
P\bigl(u^{-3/2-3\delta}\bigl|\log E\exp\bigl(\sigma
g\bigl((uI-\zz)^{-1/2}S\bigr)\bigr)\bigr|>x\bigr)\leq e^{-\lambda x^2}
\]
for $u$ sufficiently large.
\end{lemma}
\begin{pf}
Note that $g(t)$ is a mean zero Gaussian random field
with $\operatorname{Var}(g(t))=O(|t|^6)$ and $|S|\leq\kappa u^{\delta}$. A direct
application of the Borel-TIS lemma yields the result of this lemma.
\end{pf}
\begin{lemma}\label{LemPair}
For each $\kk\neq\kk'$,
\[
P\bigl(I_\sigma(\Xi_{\varepsilon
,\kk})>b,I_\sigma(\Xi_{\varepsilon,\kk'})>b\bigr)=O(1)u^{d-1}
e^{- u^2/2-\Theta(u)}.
\]
\end{lemma}
\begin{pf} 
Without loss of generality, we consider $\Xi_{\varepsilon
,\kk'}=\Xi_{\varepsilon}$. If $\Xi_{\varepsilon}$
and~$\Xi_{\varepsilon,\kk'}$ are connected to each other,
\begin{eqnarray*}
P\bigl( I_\sigma(\Xi_{\varepsilon,\kk
})>b,I_\sigma(\Xi
_{\varepsilon})>b\bigr) &\leq& P\bigl(I_\sigma(\Xi_{\varepsilon}\cup
\Xi_{\varepsilon,\kk%
})>2b\bigr) \\
&=&O(1)m(\Xi_{\varepsilon})u^{d-1}e^{-(1/2)(u+\log2+o(1))^{2}} \\
&=&O(1)u^{d-1}e^{-({1/2})u^{2}-(\log2+o(1))u}.
\end{eqnarray*}
The second step is an application of Theorem
\ref{PropSmall}. If $\Xi_{\varepsilon}$ and $\Xi
_{\varepsilon,\kk'}$ are not connected, that is, $\inf_{s\in
\Xi_\varepsilon, t\in\Xi_{\varepsilon,\kk'}}|s-t|\geq
\varepsilon= \kappa u^{-1/2 + \delta}$, then
\begin{eqnarray*}
&&P\bigl( I_\sigma(\Xi_{\varepsilon,\kk})>b,I_\sigma
(\Xi_{\varepsilon})>b\bigr) \\[-1.3pt]
&&\qquad\leq P\Bigl(\sup_{t\in\Xi
_{\varepsilon}}f(t)>u-c\log u,\sup_{t\in\Xi
_{\varepsilon,\kk}}f(t)>u-c\log u\Bigr) \\[-1.3pt]
&&\qquad\leq P\Bigl(\sup_{t\in\Xi_{\varepsilon},s\in\Xi_{\varepsilon,\kk%
}}f(t)+f(s)>2u-2c\log u\Bigr) \\[-1.3pt]
&&\qquad\leq O(1)P\biggl( Z>\frac{2u-2c\log
u +O(1)}{\sqrt{4-\Theta(1)u^{-1+2\delta
}}}\biggr) \\[-1.3pt]
&&\qquad= O(1)u^{d-1}e^{-({1/2})u^{2}-\Theta(1)u^{1+2\delta}},
\end{eqnarray*}
where $Z$ is a standard Gaussian random variable. The last inequality
is an application of the Borel-TIS lemma in Lemma
\ref{LemBorel}.
\end{pf}
\begin{lemma}\label{LemAD}
Let $\cA$ and $\cD$ be defined in (\ref{DefABD}) and $b_0$ be
defined in (\ref{B0}). Then,
\[
P( \cD>b_{0},\cA>b_{0}) =
o\bigl(P(\cA> b)\bigr).
\]
\end{lemma}
\begin{pf}
Similar to the second case in the proof of Lemma \ref{LemPair}, we
have
\begin{eqnarray*}
&&P( D>b_{0},A>b_{0}) \\[-1.3pt]
&&\qquad\leq P\Bigl( \sup_{\Xi_{\varepsilon}}f(t)>u-c_{1}\log
u,\sup_{\bigcup_{\kk\in C^{+}\setminus\{\mathbf0,\kk\}}\Xi
_{\varepsilon,\kk}}f(t)>u-c_{1}\log
u\Bigr) \\[-1.3pt]
&&\qquad\leq P\Bigl( \sup_{s\in\Xi_{\varepsilon},t\in\bigcup
_{\kk\in C^{+}\setminus\{\mathbf0,\kk\}}\Xi_{\varepsilon
,\kk}}f(s)+f(t)>2u-2c_{1}\log
u\Bigr) \\[-1.3pt]
&&\qquad\leq u^{\alpha}e^{-{(2u-2c_{1}\log u)^{2}}/({2(4-2\kappa
^{2}u^{-1+2\delta})})}\leq
e^{-{u^{2}}/{2}-\Theta(1)u^{1+2\delta}}.
\end{eqnarray*}
The conclusion follows immediately.
\end{pf}
\begin{lemma}\label{Lemlast}
Let $\cA$, $\cB$ and $\cD$ be defined in (\ref{DefABD}) and
$b_0$ be defined in (\ref{B0}). Then,
\[
P(\cA+ \cB> b-b_0 , \cA> b_0 , \cB+ \cD> b_0)=o(1) P(\cA>
b).
\]
\end{lemma}
\begin{pf}
Note that there exists $c'>0$ such that
\begin{eqnarray*}
&&P( \mathcal{A}+\mathcal{B}>b-b_{0},\mathcal
{A}>b_{0},\mathcal{B}+%
\mathcal{D}>b_{0}) \\[-1.3pt]
&&\qquad\leq P\Bigl( \mathcal{A}+\mathcal{B}>b-b_{0},\sup_{\Xi
_{\varepsilon}}f(t)>u-c^{\prime}\log u,\\[-1.3pt]
&&\qquad\hspace*{65.3pt}\sup_{\bigcup_{\kk
\in\mathcal{C}^{+}\setminus\{\mathbf{0}
\}}\Xi_{\varepsilon, \kk}}f(t)>u-c^{\prime}\log u\Bigr).
\end{eqnarray*}
It suffices to show that the right-hand side of the above
inequality is\break $o(1) P(\cA>b)$ and also $o(1) P(\cA+ \cB>b)$.
In order to do so, we will borrow part of the derivations in
the proof of Theorem \ref{PropSmall}. Let $u_*$ solve
\[
(2\pi/\sigma)^{d/2}u_*^{-d/2}e^{\sigma u_*}=b-b_0.
\]
Note that because $b_0= u^{-1-d/2} b$, we have
$|u-u_*|=o(u^{-1})$ and $e^{-u^2/2} = (1+o(1))e^{-u_*^2/2}$. By
the results in (\ref{ApxL1}), (\ref{ApxL2}) and (\ref{ApxL3}),
we have
\begin{eqnarray*}
\hspace*{-5.5pt}&&P\Bigl(
\mathcal{A}+\mathcal{B}>b-b_{0},\sup_{\Xi_{\varepsilon
}}f(t)>u-c^{\prime}\log u,\sup_{\bigcup_{\kk\in\mathcal
{C}^{+}\setminus\{%
\mathbf{0}\}}\Xi_{\varepsilon,\kk}}f(t)>u-c^{\prime}\log u\Bigr)
\\
\hspace*{-5.5pt}&&\hspace*{-3pt}\qquad=o(1)P(\mathcal{A}+\mathcal{B}>b-b_{0}) \\
\hspace*{-5.5pt}&&\hspace*{-5pt}\qquad\quad{}
+\bigl(1+o(1)\bigr)\sigma^{-1}\det(\Gamma)^{-1/2}\det(\mu
_{22})^{1/2}u^{-1}\\
\hspace*{-5.5pt}&&\hspace*{-5pt}\hspace*{11pt}\qquad\quad{}\times e^{-({1/2})u^{2}+
({1}/({8\sigma^2}))\mathbf{1}^{\top}\mu
_{22}\mathbf{1}+({1}/({8\sigma^2}))\sum_{i}\partial
_{iiii}C(0)} \\
\hspace*{-5.5pt}&&\hspace*{-5pt}\hspace*{11pt}\qquad\quad{}
\times\int_{\{|y|_{\infty}\leq3\kappa u^{1/2+\delta}\}}P\Bigl(
A^{\prime}>o_{p}(1),\sup_{\Xi_{\varepsilon
}}E(t)+g(t)>u_{\ast}-c^{\prime}\log u_{\ast},\\
\hspace*{-5.5pt}&&\hspace*{-5pt}\qquad\quad\hspace*{125.2pt}\sup_{\bigcup
_{\kk\in\mathcal{C}^{+}\setminus\{\mathbf{0} \}}\Xi
_{\varepsilon,\kk}}E(t)+g(t)>u_{\ast}-c^{\prime}\log u_{\ast
}\Bigr) \\
\hspace*{-5.5pt}&&\qquad\quad\hspace*{-5pt}\hspace*{11pt}{}
\times\exp\biggl\{-\biggl[\frac{A'}{\sigma}+\frac{B^{\top}B}{2}\\
\hspace*{-5.5pt}&&\qquad\quad\hspace*{-5pt}\hspace*{57pt}{}+\frac
{(-{A'}/ ({\sigma u})+\mu
_{20}\mu_{22}^{-1/2}B+({1}/({2\sigma}))\mu
_{20}\mathbf{1})^{2}}{2(1-\mu_{20}\mu
_{22}^{-1}\mu_{02})}\biggr]\biggr\}\,dA^{\prime}\,dB\,dy.
\end{eqnarray*}
Note that the only change in the above display from
(\ref{ApxL1}) is the probability inside the integral. In what
follows, we show that it is almost always $o(1)$. Note that
$\operatorname{Var}(g(t))=O(|t|^6)$. Therefore, for any $f(0)<u +
u^{\varepsilon_0}$ with $\varepsilon_0 < \delta/2$, if
%
%
\begin{eqnarray}\label{cond}
\sup_{\Xi_\varepsilon} E(t) &<& u -
c'\log u -\Theta(u^{-1})\quad \mbox{or }\nonumber\\[-8pt]\\[-8pt]
\sup_{\Xi_{3\varepsilon}\setminus\Xi_{\varepsilon}} E(t) &<& u -
c'\log u
-\Theta(u^{-1}),\nonumber
\end{eqnarray}
then
\[
P\Bigl( \sup_{\Xi_{\varepsilon
}}E(t)+g(t)>u-c^{\prime}\log u,\sup_{\bigcup
_{\kk\in\mathcal{C}^{+}\setminus\{\mathbf{0}\}}\Xi_{\varepsilon
,\kk}}E(t)+g(t)>u-c^{\prime}\log u\Bigr)=o(1).
\]
This fact implies that $g(t)$ can be basically ignored. Therefore, it
is useful to keep in mind that ``$E(t) \approx f(t)$.''

Since
\[
P\Bigl(\sup_T f(t)> u + u^{-1+ \varepsilon_0}\Bigr) = o(1) P(\cA
+ \cB> b-b_0),
\]
we only need to consider the case that
$\sup_T E(t) \leq u + u^{-1+ \varepsilon_0}$. Given the form
\[
E(t)=u_*-w+y^{\top}t+\tfrac{1}{2}t^{\top
}(-u_*I+\zz)t+g_{3}(t)+g_{4}(t) + R(t),
\]
which is asymptotically quadratic in $\Xi_{3\varepsilon}$, let
$t^*= \arg\sup_{\Xi_{3\varepsilon}} E(t)$. On the set that
$\sup_T E(t) \leq u + u^{-1+ \varepsilon_0}$, we have
\[
\sup_{|t-t^*|> 2 u^{-1/2+\varepsilon_0/2}} E(t) < u -c'\log u
-\Theta(u^{-1}).
\]
Let $\partial\Xi_\varepsilon$ be the
border of $\Xi_\varepsilon$. Then
%
%
\begin{eqnarray}\label{both}\sup_{\Xi_\varepsilon} E(t) &>& u -
c'\log u -\Theta(u^{-1})\quad \mbox{and}\nonumber\\[-8pt]\\[-8pt]
\sup_{\Xi_{3\varepsilon}\setminus\Xi_{\varepsilon}} E(t) &>&u -
c'\log u
-\Theta(u^{-1}),\nonumber
\end{eqnarray}
only when $\inf_{t\in\partial
\Xi_\varepsilon} |t-t^*|< u_*^{-1/2+\varepsilon_0}$. This
implies that
\[
\inf_{t\in\partial\Xi_\varepsilon}
|t-(u_*I-\zz)^{-1} y|< u_*^{-1/2+\varepsilon_0}.
\]
Therefore, $t^*=\arg\sup f(t)$ must be very closed to the boundary of
$\Xi_\varepsilon$ so as to have (\ref{both}) hold.

Therefore, for all $\varepsilon_0 <\delta$
\begin{eqnarray*}
\hspace*{-4pt}&&P\Bigl( \mathcal{A}+\mathcal{B}>b-b_{0},\sup_{\Xi
_{\varepsilon
}}f(t)>u-c^{\prime}\log u,\\
\hspace*{-4pt}&&\qquad\hspace*{2pt}\sup_{ \bigcup_{\kk\in\mathcal
{C}^{+}\setminus\{%
\mathbf{0}\}}\Xi_{\varepsilon,\kk}}E(t)+g(t)>u-c^{\prime
}\log u\Bigr)
\\
\hspace*{-4pt}&&\qquad=o(1)P(\mathcal{A}+\mathcal{B}>b-b_{0}) \\
\hspace*{-4pt}&&\hspace*{-3pt}\qquad\quad{} +\bigl(1+o(1)\bigr)\sigma^{-1}\det(\Gamma)^{-1/2}\det(\mu
_{22})^{1/2}u^{-1}\\
\hspace*{-4pt}&&\qquad\qquad\hspace*{-4pt}{}\times e^{-({1/2})u^{2}+({1}/({8\sigma^2}))
\mathbf{1}^{\top}\mu
_{22}\mathbf{1}+({1}/({8\sigma^2}))\sum_{i}\partial
_{iiii}C(0)} \\
\hspace*{-4pt}&&\qquad\qquad\hspace*{-4pt}{}\times\int_{\inf_{t\in\partial\Xi_{\varepsilon}}|t-(u_{\ast}I-\zz%
)^{-1}y|<u_{\ast}^{-1/2+\varepsilon_{0}}}\hspace*{-3pt}
P\Bigl( A^{\prime}>o_{p}(1),\\
\hspace*{-4pt}&&\hspace*{198pt}\sup_{\Xi_{\varepsilon
}}E(t)+g(t)>u-c^{\prime}\log u,\\
\hspace*{-4pt}&&\hspace*{163pt}\sup_{ \bigcup_{\kk\in\mathcal{
C}^{+}\setminus\{\mathbf{0}\}}\Xi_{\varepsilon
,\kk}}E(t)+g(t)
>u - c^{\prime}\log u\Bigr) \\
\hspace*{-4pt}&&\qquad\qquad\hspace*{-4pt}{}\times\exp\biggl\{-\biggl[\frac{A'}\sigma+\frac{B^{\top}B}{2}\\
\hspace*{-4pt}&&\qquad\quad\hspace*{54pt}{}+\frac{(-
{A'} /({\sigma u})+\mu
_{20}\mu_{22}^{-1 /2}B+({1}/({2\sigma}))\mu
_{20}\mathbf{1})^{2}}{2(1-\mu_{20}\mu
_{22}^{-1}\mu
_{02})}\biggr]\biggr\}\,dA'\,dB\,dy\\
\hspace*{-4pt}&&\qquad= o(1)P(\mathcal{A}+\mathcal{B}>b-b_{0}).
\end{eqnarray*}
The last equation is because
\begin{eqnarray*}
&&\mes\Bigl(\Bigl\{y\dvtx{\inf_{t\in\partial\Xi_{\varepsilon}}}|t-(u_{\ast}I-\zz
)^{-1}y|<u_{\ast}^{-1/2+\varepsilon_{0}}\Bigr\}\Bigr) \\
&&\qquad=
o\bigl(\mes\bigl(\{y\dvtx|(uI-\zz)^{-1/2}y|_\infty\leq\kappa
u^{\delta}-u^{\delta/2}\}\bigr)\bigr).
\end{eqnarray*}
Hereby, we conclude the proof.
\end{pf}
\end{appendix}


%

%
\printaddresses

\end{document}